\documentclass[12pt]{article}

\usepackage{amssymb,latexsym,amsmath,setspace}
\usepackage{amsthm,amsfonts}
\usepackage{graphicx}
\usepackage [left =1.3 in, top =1.0in, bottom =1.0in, right=1.35in]{geometry}

\newcommand{\norm}[1]{\left\Vert#1\right\Vert}

\newcommand{\abs}[1]{\left\vert#1\right\vert}
\newcommand{\Lip}[1]{\text {Lip}\left(#1\right)}
\newcommand{\dia}[1]{\text {diameter}\left(#1\right)}

\newcommand{\C}{\text{const}}

\newcommand{\set}[1]{\left\{#1\right\}}

\theoremstyle{plain}
\newtheorem{thm}{Theorem}
\newtheorem{cor}[thm]{Corollary}
\newtheorem{lem}[thm]{Lemma}
\newtheorem{prop}[thm]{Proposition}
\theoremstyle{definition}
\newtheorem{defn}{Definition}

\theoremstyle{remark}
\theoremstyle{plain}
\newtheorem*{main_assu}{Main Assumption}

\newcommand {\marginal} [1] {$\bigstar $\leavevmode\marginpar
    {\tiny\raggedright#1\par}}
\newcommand{\comment}[1]{}
\begin{document}

\title {\textbf{The periodic oscillation of an adiabatic piston in
two or three dimensions} \footnote {Submitted to
\emph{Communications in Mathematical Physics}}\\ preprint }

\author {Paul Wright\footnote
        {Department of Mathematics, Courant Institute of Mathematical Sciences, New York University,
        251 Mercer St.,
        New York, NY 10012 USA.  E-mail address:~\texttt{paulrite@cims.nyu.edu}}}

\date {December 2006}

\maketitle

\begin {abstract}

We study a heavy piston of mass $M$ that separates finitely many
ideal, unit mass gas particles moving in two or three dimensions.
Neishtadt and Sinai previously determined a method for finding this
system's averaged equation and showed that its solutions oscillate
periodically.  Using averaging techniques, we prove that the actual
motions of the piston converge in probability to the predicted
averaged behavior on the time scale $M^ {1/2} $ when $M$ tends to
infinity while the total energy of the system is bounded and the
number of gas particles is fixed.

\end {abstract}




\section{Introduction}\label{sct:intro}

Consider the following simple model of an adiabatic piston
separating two gas containers:  A massive piston of mass $M\gg 1$
divides a container in $\mathbb{R}^d$ into two halves. The piston
has no internal degrees of freedom and can only move along one axis
of the container. On either side of the piston there are a finite
number of ideal, unit mass, point gas particles that interact with
the walls of the container and with the piston via elastic
collisions. When $M=\infty $, the piston remains fixed in place, and
each gas particle performs billiard motion at a constant energy in
its sub-container. We make an ergodicity assumption on the behavior
of the gas particles when the piston is fixed. Then we study the
motions of the piston when the number of gas particles is fixed, the
total energy of the system is bounded, but $M$ is very large.

Heuristically, after some time, one expects the system to approach a
steady state, where the energy of the system is equidistributed
amongst the particles and the piston. However, even if we could show
that the full system is ergodic, an abstract ergodic theorem says
nothing about the time scale required to reach such a steady state.
Because the piston will move much slower than a typical gas
particle, it is natural to try to determine the intermediate
behavior of the piston by averaging techniques. By averaging over
the motion of the gas particles on a time scale chosen short enough
that the piston is nearly fixed, but long enough that the ergodic
behavior of individual gas particles is observable, we will show
that the system does not approach the expected steady state on the
time scale $M^ {1/2} $. Instead, the piston oscillates periodically,
and there is no net energy transfer between the gas particles.

This paper follows earlier work by Neishtadt and Sinai~\cite{NS04,
Sin99}. They determined that for a wide variety of Hamiltonians for
the gas particles, the averaged behavior of the piston is periodic
oscillation, with the piston moving inside an effective potential
well whose shape depends on the initial position of the piston and
the gas particles' Hamiltonians.  They pointed out that an averaging
theorem due to Anosov~\cite{Ano60,LM88}, proved for smooth systems,
should extend to this case. This paper proves that Anosov's theorem
extends to the particular gas particle Hamiltonian described above.
Thus, if we examine the actual motions of the piston with respect to
the slow time $\tau=t/M^ {1/2}$,  then, as $M\rightarrow\infty $, in
probability (with respect to Liouville measure) most initial
conditions give rise to orbits whose actual motion is accurately
described by the averaged behavior for $0\leq\tau\leq 1$, i.e.~for
$0\leq t\leq M^ {1/2}$.  Gorelyshev and
Neishtadt~\cite{GorNeishtadt06} and we~\cite{Wri06} have already
proved that when $d=1$, i.e.~when the gas particles move on a line,
the convergence of the actual motions to the averaged behavior is
uniform over all initial conditions, with the size of the deviations
being no larger than $\mathcal{O} (M^ {-1/2}) $ on the time scale
$M^ {-1/2} $.

The system under consideration in this paper is a simple model of an
adiabatic piston. The general adiabatic piston problem~\cite {Ca63},
well-known from physics, consists of the following: An insulating
piston separates two gas containers, and initially the piston is
fixed in place, and the gas in each container is in a separate
thermal equilibrium. At some time, the piston is no longer
externally constrained and is free to move. One hopes to show that
eventually the system will come to a full thermal equilibrium, where
each gas has the same pressure and temperature. Whether the system
will evolve to thermal equilibrium and the interim behavior of the
piston are mechanical problems, not adequately described by
thermodynamics~\cite{Gru99}, that have recently generated much
interest within the physics and mathematics communities.  One
expects that the system will evolve in at least two stages. First,
the system relaxes toward a mechanical equilibrium, where the
pressures on either side of the piston are equal.  In the second,
much longer, stage, the piston drifts stochastically in the
direction of the hotter gas, and the temperatures of the gases
equilibrate.  See for example~\cite{GPL03,CL02} and the references
therein. So far, rigorous results have been limited mainly to models
where the effects of gas particles recolliding with the piston can
be neglected, either by restricting to extremely short time
scales~\cite{CLS02,CLS02b} or to infinite gas
containers~\cite{Che05}.

A recent study involving some similar ideas by Chernov and
Dolgopyat~\cite{CD06} considered the motion inside a two-dimensional
domain of a single heavy, large gas particle (a disk) of mass $M\gg
1$ and a single unit mass point particle.  They assumed that for
each fixed location of the heavy particle, the light particle moves
inside a dispersing (Sinai) billiard domain. By averaging over the
strongly hyperbolic motions of the light particle, they showed that
under an appropriate scaling of space and time the limiting process
of the heavy particle's velocity is a (time-inhomogeneous) Brownian
motion on a time scale $\mathcal{O} (M^ {1/2}) $. It is not clear
whether a similar result holds for the piston problem, even for gas
containers with good hyperbolic properties, such as the Bunimovich
stadium.  In such a container the motion of a gas particle when the
piston is fixed is only nonuniformly hyperbolic because it can
experience many collisions with the flat walls of the container
immediately preceding and following a collision with the piston.

The present work provides a weak law of large numbers, and it is an
open problem to describe the sizes of the deviations for the piston
problem~\cite{CD06b}.  Although our result does not yield concrete
information on the sizes of the deviations, it is general in that it
imposes very few conditions on the shape of the gas container. Most
studies of billiard systems impose strict conditions on the shape of
the boundary, generally involving the sign of the curvature and how
the corners are put together.  The proofs in this work require no
such restrictions.  In particular, the gas container can have cusps
as corners and need satisfy no hyperbolicity conditions.

We begin in Section~\ref{sct:heuristic} by giving a physical
description of our results.  Precise assumptions and our main
result, Theorem~\ref{thm:dDpiston}, are stated in
Section~\ref{sct:main_result}, and a proof is presented in the
following sections.

\section{Physical motivation for the results}\label{sct:heuristic}

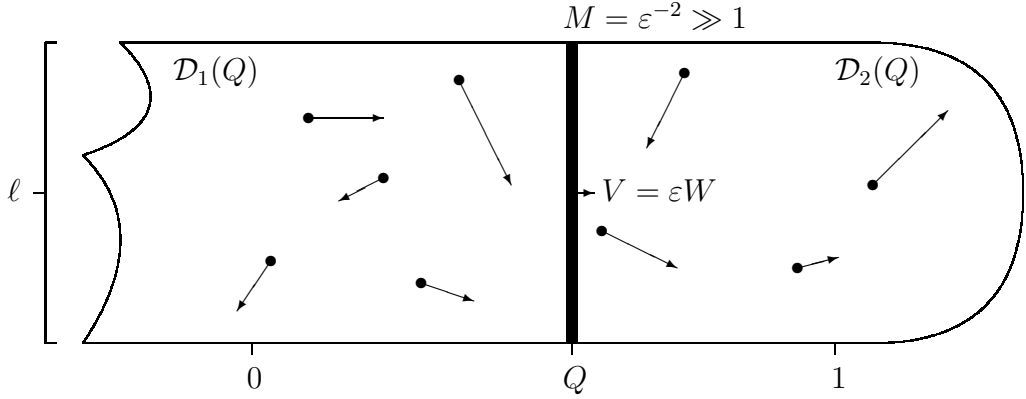
\begin{figure}
    \begin {center}
    \setlength{\unitlength}{1 cm}
    \begin{picture}(15,6)
        \put(2.5,1){\line(1,0){10}}
        \put(2.5,5){\line(1,0){10}}
        \qbezier(2,3.5)(3.5,4)(2.5,5)
        \qbezier(2,3.5)(3,2.5)(2,1)
        \put(2,1){\line(1,0){1}}
        \qbezier(12.5,1)(14.5,1)(14.5,3)
        \qbezier(14.5,3)(14.5,5)(12.5,5)
        \put(3.2,4.5){$\mathcal{D}_1(Q)$}
        \put(12,4.5){$\mathcal{D}_2(Q)$}
        \put(1.5,1){\line(0,1){4}}
        \put(1.5,1){\line(1,0){0.15}}
        \put(1.5,5){\line(1,0){0.15}}
        \put(1.5,3){\line(-1,0){0.15}}
        \put(1,2.9){$\ell$}
        \linethickness {0.15cm}
        \put(8.5,1){\line(0,1){4}}
        \thinlines
        \put(8.5,3){\vector(1,0){.3}}
        \put(8.9,2.9){$V=\varepsilon W$}
        \put(8.5,1){\line(0,-1){0.15}}
        \put(8.38,0.4){$Q$}
        \put(8.38,5.2){$M=\varepsilon ^ {-2}\gg 1$}
        \put(4.25,1){\line(0,-1){0.15}}
        \put(4.18,0.4){$0$}
        \put(12,1){\line(0,-1){0.15}}
        \put(11.95,0.4){$1$}
        \put(5,4){\circle*{.15}}
        \put(5,4){\vector(1,0){1}}
        \put(7,4.5){\circle*{.15}}
        \put(7,4.5){\vector(1,-2){0.7}}
        \put(6,3.2){\circle*{.15}}
        \put(6,3.2){\vector(-2,-1){0.6}}
        \put(6.5,1.8){\circle*{.15}}
        \put(6.5,1.8){\vector(3,-1){0.7}}
        \put(4.5,2.1){\circle*{.15}}
        \put(4.5,2.1){\vector(-2,-3){0.45}}
        \put(8.9,2.5){\circle*{.15}}
        \put(8.9,2.5){\vector(2,-1){1}}
        \put(11.5,2.0){\circle*{.15}}
        \put(11.5,2.0){\vector(4,1){0.55}}
        \put(10,4.6){\circle*{.15}}
        \put(10,4.6){\vector(-1,-2){0.5}}
        \put(12.5,3.1){\circle*{.15}}
        \put(12.5,3.1){\vector(1,1){1}}
    \end{picture}
    \end {center}
    \caption{A gas container $\mathcal{D}\subset \mathbb{R}^2 $ separated by a piston.}
    \label{fig:domain1}
\end{figure}

Before precisely stating our assumptions and results, we briefly
review the physical motivations for our results and introduce some
notation.

Consider a massive, insulating piston of mass $M$ that separates a
gas container $\mathcal{D} $ in $\mathbb{R}^d$, $ d= 2\text { or
}3$. See Figure~\ref{fig:domain1}. Denote the location of the piston
by $Q$, its velocity by $dQ/dt=V$, and its cross-sectional length
(when $ d=2$, or area, when $ d=3$) by $\ell$. If $Q$ is fixed, then
the piston divides $\mathcal{D} $ into two subdomains,
$\mathcal{D}_1(Q) =\mathcal{D}_1 $ on the left and $\mathcal{D}_2(Q)
=\mathcal{D}_2  $ on the right. By $E_i$ we denote the total energy
of the gas inside $\mathcal{D}_i$, and by $\abs{\mathcal{D}_i} $ we
denote the area (when $ d=2$, or volume, when $ d=3$) of
$\mathcal{D}_i$.

We are interested in the dynamics of the piston when the system's
total energy is bounded and $M\rightarrow \infty $.  When
$M=\infty$, the piston remains fixed in place, and each energy $E_i$
remains constant. When $M$ is large but finite, $MV^2/2$ is bounded,
and so $V=\mathcal{O} (M^ {-1/2}) $.  It is natural to define
\[
    \varepsilon=M^ {-1/2},\quad W=\frac {V} {\varepsilon},
\]
so that $W$ is of order $1$ as $\varepsilon\rightarrow 0$.  This is
equivalent to scaling time by $\varepsilon$.

If we let $P_i$ denote the pressure of the gas inside
$\mathcal{D}_i$, then heuristically the dynamics of the piston
should be governed by the following differential equation:
\[
    \frac {dQ} {dt}=V,\quad M\frac {dV} {dt}=P_1\ell-P_2\ell,
\]
i.e.
\begin {equation}\label{eq:piston_force}
    \frac {dQ} {dt}=\varepsilon W,\quad
    \frac {dW} {dt}=\varepsilon P_1\ell-\varepsilon P_2\ell.
\end {equation}
To find differential equations for the energies of the gases,
note that in a short amount of time $dt$, the change in energy
should come entirely from the work done on a gas, i.e.~the force
applied to the gas times the distance the piston has moved, because
the piston is adiabatic. Thus, one expects that
\begin {equation}\label{eq:work}
    \frac {dE_1} {dt}=-\varepsilon WP_1\ell,
    \quad
    \frac {dE_2} {dt} = +\varepsilon W P_2\ell.
\end {equation}
To obtain a closed system of differential equations, it is necessary
to insert an expression for the pressures.  Because the pressure of
an ideal gas in $d$ dimensions is proportional to the energy
density, with the constant of proportionality $2/d$, we choose to
insert
\[
    P_i=\frac {2E_i}{d\abs{\mathcal{D}_i}} .
\]
Later, we will make assumptions to justify this substitution.
However, if we accept this definition of the pressure, and define
the slow time
\[
    \tau=\varepsilon t,
\]
we obtain the following ordinary differential equations for the
four macroscopic variables of the system:
\begin {equation}\label{eq:heuristic_avg_eq}
    \frac{d}{d\tau}
    \begin {bmatrix}
    Q\\
    W\\
    E_1\\
    E_2\\
    \end {bmatrix}
    =
    \begin {bmatrix}
    W\\
    \frac{2E_1\ell}{d\abs{\mathcal{D}_1(Q)}}
    -\frac{2E_2\ell}{d\abs{\mathcal{D}_2(Q)}}\\
    -\frac{2WE_1\ell}{d\abs{\mathcal{D}_1(Q)}}\\
    +\frac{2WE_2\ell}{d\abs{\mathcal{D}_2(Q)}}\\
    \end {bmatrix}.
\end {equation}

Neishtadt and Sinai~\cite{Sin99, NS04} pointed out that the
solutions of Equation~\eqref{eq:heuristic_avg_eq} have the piston
moving according to an effective Hamiltonian.  This can be seen as
follows.  Since
\[
    \frac {\partial\abs{\mathcal{D}_1(Q)}} {\partial Q}=\ell
    =-\frac {\partial\abs{\mathcal{D}_2(Q)}} {\partial Q},
\]
$d\ln(E_i)/d\tau=-(2/d)d\ln(\abs{\mathcal{D}_i(Q)})/d\tau$, and so
\[
    E_i(\tau)=E_i(0)\left (\frac{\abs{\mathcal{D}_i(Q(0))}}
        {\abs{\mathcal{D}_i(Q(\tau))}}\right) ^ {2/d}.
\]
Hence
\[
    \frac {d^2Q} {d\tau^2}=
    \frac {2\ell} {d}
    \frac{E_1(0)\abs{\mathcal{D}_1(Q(0))}^{2/d}}{\abs{\mathcal{D}_1(Q(\tau))}^{1+2/d}}
    -\frac {2\ell} {d}
    \frac{E_2(0)\abs{\mathcal{D}_2(Q(0))}^{2/d}}{\abs{\mathcal{D}_2(Q(\tau))}^{1+2/d}},
\]
and so $(Q, W) $ behave as if they were the coordinates of a
Hamiltonian system describing a particle undergoing motion inside a
potential well. The effective Hamiltonian may be expressed as
\begin {equation}\label{eq:d_dpot}
    \frac {1} {2}W^2+
    \frac{E_1(0)\abs{\mathcal{D}_1(Q(0))}^{2/d}}
    {\abs{\mathcal{D}_1(Q)}^{2/d}}+
    \frac{E_2(0)\abs{\mathcal{D}_2(Q(0))}^{2/d}}
    {\abs{\mathcal{D}_2(Q)}^{2/d}}.
\end {equation}

The question is, do the solutions of
Equation~\eqref{eq:heuristic_avg_eq} give an accurate description of
the actual motions of the macroscopic variables when $M$ tends to
infinity?  The main result of this paper is that, for an
appropriately defined system, the answer to this question is
affirmative for $0\leq t\leq M^ {1/2}$, at least for most initial
conditions of the microscopic variables. Observe that one should not
expect the description to be accurate on time scales much longer
than $\mathcal{O} (M^ {1/2}) =\mathcal{O} (\varepsilon^ {-1})$.  The
reason for this is that, presumably, there are corrections of size
$\mathcal{O} (\varepsilon^ {2})$ in
Equations~\eqref{eq:piston_force} and \eqref{eq:work} that we are
neglecting.  On the time scale $\varepsilon^ {-1} $, these errors
roughly add up to no more than size $\mathcal{O} (\varepsilon^ {-1}
\cdot \varepsilon^ {2}=\varepsilon) $, but on a longer time scale
they should become significant.  Such higher order corrections for
the adiabatic piston were studied by Crosignani \emph{et
al.}~\cite{CD96}.

\section{Statement of the main result}\label{sct:main_result}

\subsection{Description of the model}\label{sct:model}

We begin by describing the gas container.  It is a compact,
connected billiard domain $\mathcal{D} \subset\mathbb{R}^d$ with a
piecewise $\mathcal{C} ^3$ boundary, i.e.~$\partial\mathcal{D} $
consists of a finite number of $\mathcal{C} ^3$ embedded
hypersurfaces, possibly with boundary and a finite number of corner
points.  The container consists of a ``tube,'' whose perpendicular
cross-section $\mathcal{P} $ is the shape of the piston, connecting
two disjoint regions.  $\mathcal{P} \subset \mathbb{R} ^ {d-1} $ is
a compact, connected domain whose boundary is piecewise $\mathcal{C}
^3$.  Then the ``tube'' is the region $[0,1]\times
\mathcal{P}\subset\mathcal{D} $ swept out by the piston for $ 0\leq
Q\leq 1$, and $[0,1]\times
\partial\mathcal{P}\subset\partial\mathcal{D} $.  If $d=2$, $\mathcal{P} $
is just a closed line segment, and the ``tube'' is a rectangle.  If
$ d=3$, $\mathcal{P} $ could be a circle, a square, a pentagon, etc.

Our fundamental assumption is as follows:
\begin{main_assu}
For almost every $Q\in [0,1]$ the billiard flow of a single particle
on an energy surface in either of the two subdomains
$\mathcal{D}_i(Q)$ is ergodic (with respect to the invariant
Liouville measure).
\end {main_assu}
\noindent If $ d=2$, the domain could be the Bunimovich
stadium~\cite{Bun79}. Another possible domain is indicated in Figure
\ref{fig:domain1}.  Polygonal domains satisfying our assumptions can
also be constructed~\cite{Vorobets_1997}. Suitable domains in $d=3$
dimensions can be constructed using a rectangular box with shallow
spherical caps adjoined~\cite{BunimovichRehacek1998}.  Note that we
make no assumptions regarding the hyperbolicity of the billiard flow
in the domain.

The Hamiltonian system we consider consists of the massive piston of
mass $M$ located at position $Q$, as well as $ n_1+ n_2 $ gas
particles, $ n_1$ in $\mathcal{D}_1$ and $n_2$ in $\mathcal{D}_2$.
Here $n_1$ and $n_2$ are fixed positive integers. For convenience,
the gas particles all have unit mass, though all that is important
is that each gas particle has a fixed mass. We denote the positions
of the gas particles in $\mathcal{D}_i$ by $q_{ i,j}$, $1\leq j\leq
n_i$.  The gas particles are ideal point particles that interact
with $\partial\mathcal{D} $ and the piston by hard core, elastic
collisions. Although it has no effect on the dynamics we consider,
for convenience we complete our description of the Hamiltonian
dynamics by specifying that the piston makes elastic collisions with
walls located at $Q=0,\: 1$ that are only visible to the piston.  We
denote velocities by $dQ/dt=V=\varepsilon W$ and $dq_{ i,j}/dt=v_{
i,j}$, and we set
\[
    E_{ i,j}=v_{ i,j}^2/2,\qquad E_i=\sum_{ j=1} ^
    {n_i} E_{ i,j}.
\]
Our system has $d(n_1+n_2)+1$ degrees of
freedom, and so its phase space is $(2d(n_1+n_2)+2)$-dimensional.

We let
\[
    h(z)=h=(Q,W,E_{1,1},E_{1,2},\cdots,E_{1,n_1},E_{2,1},E_{2,2},\cdots,E_{2,n_2}),
\]
so that $h$ is a function from our phase space to
$\mathbb{R}^{n_1+n_2+2}$.  We often abbreviate
$h=(Q,W,E_{1,j},E_{2,j})$, and we refer to $h$ as consisting of the
slow variables because these quantities are conserved when
$\varepsilon=0$.  We let $h_\varepsilon(t,z)=h_\varepsilon(t) $
denote the actual motions of these variables in time for a fixed
value of $\varepsilon$. Here $z$ represents the initial condition in
phase space, which we usually suppress in our notation.  One should
think of $h_\varepsilon(\cdot) $ as being a random variable that
takes initial conditions in phase space to paths (depending on the
parameter t) in $\mathbb{R}^{n_1+n_2+2}$.

\subsection {The averaged equation}

From the work of Neishtadt and Sinai~\cite{NS04}, one can derive
\begin {equation}\label{eq:d_davg}
    \frac{d}{d\tau}
    \begin {bmatrix}
    Q\\
    W\\
    E_{1,j}\\
    E_{2,j}\\
    \end {bmatrix}
    =\bar H(h):=
    \begin {bmatrix}
    W\\
    \frac{2E_1\ell}{d\abs{\mathcal{D}_1(Q)}}
    -\frac{2E_2\ell}{d\abs{\mathcal{D}_2(Q)}}\\
    -\frac{2WE_{1,j}\ell}{d\abs{\mathcal{D}_1(Q)}}\\
    +\frac{2WE_{2,j}\ell}{d\abs{\mathcal{D}_2(Q)}}\\
    \end {bmatrix}
\end {equation}
as the averaged equation (with respect to the slow time
$\tau=\varepsilon t$) for the slow variables. Later, in Section
\ref{sct:heuristic2}, we will give another heuristic derivation of
the averaged equation that is more suggestive of our proof.  As in
Section~\ref{sct:heuristic}, the solutions of
Equation~\eqref{eq:d_davg} have $(Q, W) $ behaving as if they were
the coordinates of a Hamiltonian system describing a particle
undergoing motion inside a potential well. The effective Hamiltonian
is given by Equation~\eqref{eq:d_dpot}.

Let $\bar{h} (\tau,z)=\bar{h} (\tau) $ be the solution of
\[
\frac {d\bar{h}}{d\tau} =\bar {H} (\bar {h}),\qquad \bar {h} (0)
=h_\varepsilon(0).
\]
Again, think of $\bar h(\cdot) $ as being a random variable.

\subsection{The main result}

The solutions of the averaged equation approximate the motions of
the slow variables, $h_\varepsilon(t) $, on a time scale
$\mathcal{O} (1/\varepsilon) $ as $\varepsilon\rightarrow 0$.
Precisely, fix a compact set $\mathcal{V}\subset \mathbb
R^{n_1+n_2+2}$ such that $h\in \mathcal{V} \Rightarrow
Q\subset\subset (0,1),W\subset\subset \mathbb R$, and $E_{i,j}
\subset\subset (0,\infty)$ for each $i$ and $j$.\footnote { We have
introduced this notation for convenience.  For example, $h\in
\mathcal{V} \Rightarrow Q\subset\subset (0,1) $ means that there
exists a compact set $A \subset (0,1) $ such that $h\in \mathcal{V}
\Rightarrow Q\in A $, and similarly for the other variables.}  We
will be mostly concerned with the dynamics when $h\in\mathcal{V} $.
Define
\[
\begin {split}
    Q_{min}&=\inf_{h\in\mathcal{V}}Q,\qquad
    Q_{max}=\sup_{h\in\mathcal{V}}Q,
    \\
    E_{min}&=\inf_{h\in\mathcal{V}}\frac{1}{2}W^2+E_1+E_2,\qquad
    E_{max}=\sup_{h\in\mathcal{V}}\frac{1}{2}W^2+E_1+E_2.
\end {split}
\]
For a fixed value of $\varepsilon >0$, we only consider the dynamics
on the invariant subset of phase space defined by
\[
\begin {split}
    \mathcal{M}_\varepsilon =
    \{(Q,V,q_{i,j},v_{i,j})\in\mathbb{R}^ {2d(n_1+n_2)+2}:
    Q\in [0,1],\;q_{i,j}\in\mathcal{D}_{i}(Q),
    &
    \\
    E_{min}\leq \frac{M}{2}V^2+E_1+E_2\leq E_{ max}\}
    &.
\end {split}
\]
Let $P_\varepsilon$ denote the probability measure obtained by
restricting the invariant Liouville measure to
$\mathcal{M}_\varepsilon $.  Define the stopping time
\[
    T_\varepsilon(z) =T_\varepsilon =\inf \{\tau\geq 0: \bar {h}
    (\tau)\notin \mathcal{V} \text { or } h_\varepsilon(\tau
    /\varepsilon) \notin \mathcal{V} \}.
\]

\begin {thm}\label{thm:dDpiston}

If $\mathcal{D} $ is a gas container in $d=2$ or $3$ dimensions
satisfying the assumptions in Subsection~\ref{sct:model} above, then
for each $T>0$,
\[
    \sup_{0\leq\tau\leq T\wedge
    T_\varepsilon}\abs{h_\varepsilon(\tau/\varepsilon)-\bar{h}(\tau)}
    \rightarrow 0 \text { in probability as } \varepsilon=M^
    {-1/2}\rightarrow 0,
\]
i.e.~for each fixed $\delta>0$,
\[
    P_\varepsilon\left(\sup_{0\leq\tau\leq T\wedge
    T_\varepsilon}\abs{h_\varepsilon(\tau/\varepsilon)-\bar{h}(\tau)}\geq\delta
    \right)
    \rightarrow 0\text { as } \varepsilon=M^
    {-1/2}\rightarrow 0.
\]
\end{thm}

It should be noted that the stopping time in the above result is not
unduly restrictive.  If the initial pressures of the two gasses are
not too mismatched, then the solution to the averaged equation is a
periodic orbit, with the effective potential well keeping the piston
away from the walls.  Thus, if the actual motions follow the
averaged solution closely for $0\leq\tau\leq T\wedge T_\varepsilon$,
and the averaged solution stays in $\mathcal{V} $, it follows that
$T_\varepsilon
>T$.

The techniques of this paper should immediately generalize to prove
the analogue of Theorem~\ref{thm:dDpiston} above in the nonphysical
dimensions $d>3$, although we do not pursue this here.

\section{Preparatory material concerning a two-dimensional gas
    container with only one gas particle on each side}
\label{sct:2dprep}

\begin{figure}
    \begin {center}
    \setlength{\unitlength}{0.7 cm}
    \begin{picture}(15,6)
        \put(2.5,1){\line(1,0){10}}
        \put(2.5,5){\line(1,0){10}}
        \qbezier(2,3.5)(3.5,4)(2.5,5)
        \qbezier(2,3.5)(3,2.5)(2,1)
        \put(2,1){\line(1,0){1}}
        \qbezier(12.5,1)(14.5,1)(14.5,3)
        \qbezier(14.5,3)(14.5,5)(12.5,5)
        \put(3.2,4.5){$\mathcal{D}_1$}
        \put(12,4.5){$\mathcal{D}_2$}
        \put(13.3,5){\vector(0,1){0.5}}
        \put(13.23,5.7){$\parallel$}
        \put(13.3,5){\vector(1,0){0.5}}
        \put(13.95,4.9){$\perp$}
        \put(1.5,1){\line(0,1){4}}
        \put(1.5,1){\line(1,0){0.15}}
        \put(1.5,5){\line(1,0){0.15}}
        \put(1.5,3){\line(-1,0){0.15}}
        \put(1,2.9){$\ell$}
        \linethickness {0.1cm}
        \put(8.5,1){\line(0,1){4}}
        \thinlines
        \put(8.5,3){\vector(1,0){.4}}
        \put(9.0,2.9){$V=\varepsilon W$}
        \put(8.5,1){\line(0,-1){0.15}}
        \put(8.38,0.3){$Q$}
        \put(8.38,5.2){$M=\varepsilon ^ {-2}\gg 1$}
        \put(4.25,1){\line(0,-1){0.15}}
        \put(4.15,0.3){$0$}
        \put(12,1){\line(0,-1){0.15}}
        \put(11.95,0.3){$1$}
        \put(5,4){\circle*{.15}}
        \put(5,4){\vector(1,2){.4}}
        \put(5.05,3.65){$q_1$}
        \put(5.45,4.5){$v_1$}
        \put(11.5,2.5){\circle*{.15}}
        \put(11.5,2.5){\vector(-1,-3){.3}}
        \put(11.65,2.3){$q_2$}
        \put(11.35,1.5){$v_2$}
    \end{picture}
    \end {center}
    \caption{A choice of coordinates on phase space.}
    \label{fig:domain}
\end{figure}
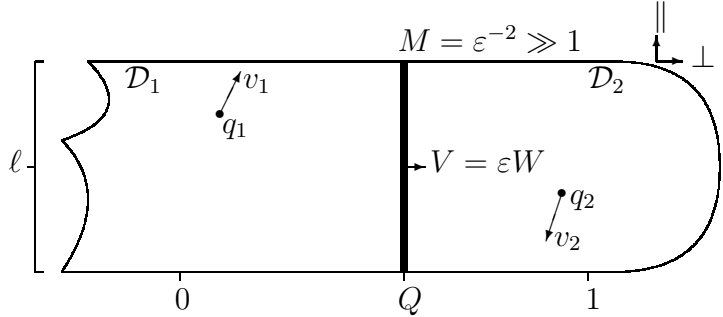

Our results and techniques of proof are essentially independent of
the dimension and the fixed number of gas particles on either side
of the piston. Thus, we focus on the case when $d=2$ and there is
only one gas particle on either side.  Later, in Section
\ref{sct:generalization}, we will indicate the simple modifications
that generalize our proof to the general situation.  For clarity, in
this section and next, we denote $ q_{ 1,1} $ by $ q_1$, $v_{2, 1} $
by $ v_2$, etc. We decompose the gas particle coordinates according
to whether they are perpendicular to or parallel to the piston's
face, for example $q_1= (q_1^\perp,q_1^\parallel)$. See Figure
\ref{fig:domain}.

The Hamiltonian dynamics define a flow on our phase space. We denote
this flow by $z_\varepsilon(t,z) =z_\varepsilon(t)$, where
$z=z_\varepsilon(0,z) $. One should think of $z_\varepsilon(\cdot) $
as being a random variable that takes initial conditions in phase
space to paths in phase space.  Then $h_\varepsilon(t)
=h(z_\varepsilon(t)) $.  By the change of coordinates
$W=V/\varepsilon$, we may identify all of the
$\mathcal{M}_\varepsilon$ defined in Section~\ref{sct:main_result}
with the space
\[
\begin {split}
    \mathcal{M} =
    \{(Q,W,q_1,v_1,q_2,v_2)\in\mathbb{R}^ {10}:
    Q\in [0,1],\;q_1\in\mathcal{D}_1(Q),\;
    q_2\in\mathcal{D}_2(Q),\;
    &
    \\
    E_{min}\leq \frac{1}{2}W^2+E_1+E_2\leq E_{ max}\}
    &.
\end {split}
\]
and all of the $P_\varepsilon$ with the probability measure $P$ on
$\mathcal{M}$, which has the density
\[
    dP=\C \,dQdWdq_1^\perp dq_1^\parallel dv_1^\perp dv_1^\parallel
    dq_2^\perp dq_2^\parallel dv_2^\perp dv_2^\parallel.
\]
(Throughout this work we will use $\C$ to represent generic
constants that are independent of $\varepsilon$.) We will assume
that these identifications have been made, so that we may consider
$z_\varepsilon(\cdot) $ as a family of measure preserving flows on
the same space that all preserve the same probability measure.  We
denote the components of $z_\varepsilon(t) $ by $Q_\varepsilon(t) $,
$q_{1,\varepsilon}^\perp(t) $, etc.

The set $\{z\in\mathcal{M}:q_1=Q=q_2\}$ has co-dimension two, and so
$\bigcup_t z_\varepsilon(t)\{q_1=Q=q_2\}$ has co-dimension one,
which shows that only a measure zero set of initial conditions will
give rise to three particle collisions.  We ignore this and other
measures zero events, such as gas particles hitting singularities of
the billiard flow, in what follows.

Now we present some background material, as well as some lemmas that
will assist us in our proof of Theorem \ref{thm:dDpiston}.  We begin
by studying the billiard flow of a gas particle when the piston is
infinitely massive.  Next we examine collisions between the gas
particles and the piston when the piston has a large, but finite,
mass.  Then we present a heuristic derivation of the averaged
equation that is suggestive of our proof.  Finally we prove a lemma
that allows us to disregard the possibility that a gas particle will
move nearly parallel to the piston's face -- a situation that is
clearly bad for having the motions of the piston follow the
solutions of the averaged equation.

\subsection{Billiard flows and maps in two dimensions}\label{sct:billiard}

In this section, we study the billiard flows of the gas particles
when $M=\infty $ and the slow variables are held fixed at a specific
value $h\in\mathcal{V} $.  We will only study the motions of the
left gas particle, as similar definitions and results hold for the
motions of the right gas particle.  Thus we wish to study the
billiard flow of a point particle moving inside the domain
$\mathcal{D}_1$ at a constant speed $\sqrt{2E_1} $.  The results of
this section that are stated without proof can be found
in~\cite{CM06}.

Let $\mathcal{TD}_1$ denote the tangent bundle to $\mathcal{D}_1$.
The billiard flow takes place in the three-dimensional space
$\mathcal{M}_h^1=\mathcal{M}^1=\{(q_1,v_1)\in\mathcal{TD}_1:q_1\in\mathcal{
D}_1,\; \abs{v_1}=\sqrt{2E_1}\}/\sim$. Here the quotient means that
when $q_1\in\partial\mathcal{ D}_1$, we identify velocity vectors
pointing outside of $\mathcal{D}_1$ with those pointing inside
$\mathcal{D}_1$ by reflecting through the tangent line to
$\partial\mathcal{D}_1$ at $ q_1$, so that the angle of incidence
with the unit normal vector to $\partial\mathcal{ D}_1$ equals the
angle of reflection. Note that most of the quantities defined in
this subsection depend on the fixed value of $h$.  We will usually
suppress this dependence, although, when necessary, we will indicate
it by a subscript $h$.  We denote the resulting flow by
$y(t,y)=y(t)$, where $y(0,y)=y$. As the billiard flow comes from a
Hamiltonian system, it preserves Liouville measure restricted to the
energy surface.  We denote the resulting probability measure by $\mu
$. This measure has the density
$d\mu=dq_1dv_1/(2\pi\sqrt{2E_1}\abs{\mathcal{D}_1} ) $. Here $dq_1$
represents area on $\mathbb{R}^2$, and $ dv_1$ represents length on
$S^1_{\sqrt{2E_1}}=\set{v_1\in\mathbb{R}^2:\abs{v_1}=\sqrt{2E_1}}$.

There is a standard cross-section to the billiard flow, the
collision cross-section
$\Omega=\{(q_1,v_1)\in\mathcal{TD}_1:q_1\in\partial\mathcal{ D}_1,\;
\abs{v_1}=\sqrt{2E_1}\}/\sim$.  It is customary to parameterize
$\Omega$ by $\{x=(r,\varphi):r\in\partial\mathcal{D}_1,\:\varphi\in
[-\pi /2,+\pi /2]\}$, where $r$ is arc length and $\varphi$
represents the angle between the outgoing velocity vector and the
inward pointing normal vector to $\partial\mathcal{D}_1$. It follows
that $\Omega$ may be realized as the disjoint union of a finite
number of rectangles and cylinders.  The cylinders correspond to
fixed scatterers with smooth boundary placed inside the gas
container.
 If $F:\Omega\circlearrowleft$ is the collision map, i.e.~the
return map to the collision cross-section, then $F$ preserves the
projected probability measure $\nu $, which has the density
$d\nu=\cos\varphi\, d\varphi \, dr/(2\abs{\partial\mathcal{D}_1}) $.
Here $\abs{\partial\mathcal{D}_1}$ is the length of
$\partial\mathcal{D}_1$.

We suppose that the flow is ergodic, and so $F$ is an invertible,
ergodic measure preserving transformation. Because
$\partial\mathcal{D}_1$ is piecewise $\mathcal{C} ^3$, $F$ is
piecewise $\mathcal{C} ^2$, although it does have discontinuities
and unbounded derivatives near discontinuities corresponding to
grazing collisions.  Because of our assumptions on $\mathcal{D}_1$,
the free flight times and the curvature of $\partial\mathcal{D}_1$
are uniformly bounded.  It follows that if $x\notin \partial\Omega
\cup F^{-1} (\partial\Omega) $, then $F$ is differentiable at $x$,
and
\begin {equation}\label{eq:2d_derivative_bound}
    \norm {DF(x)}\leq\frac {\C} {\cos \varphi(Fx)},
\end{equation}
where $\varphi(Fx) $ is the value of the $\varphi$ coordinate at the
image of $x$.

Following the ideas in Appendix \ref{sct:inducing}, we induce $F$ on
the subspace $\hat\Omega$ of $\Omega$ corresponding to collisions
with the (immobile) piston.  We denote the induced map by $\hat F$
and the induced measure by $\hat \nu$. We parameterize $\hat\Omega$
by $\{(r,\varphi):0\leq r\leq \ell,\:\varphi\in [-\pi/2,+\pi/2]\}$.
As $\nu \hat\Omega=\ell/\abs{\partial\mathcal{D}_1} $, it follows
that $\hat\nu $ has the density $d\hat \nu=\cos\varphi\, d\varphi \,
dr/(2\ell) $.

For $x\in\Omega $, define $\zeta x $ to be the free flight time,
i.e.~the time it takes the billiard particle traveling at speed
$\sqrt{2E_1} $ to travel from $x$ to $Fx$. If $x\notin
\partial\Omega \cup F^{-1} (\partial\Omega) $,
\begin {equation}
\label{eq:2d_time_derivative}
    \norm {D\zeta (x)}\leq\frac {\C} {\cos \varphi(Fx)}.
\end {equation}
Santal\'{o}'s formula~\cite{San76,Chernov1997} tells us that
\begin {equation}
\label{eq:2d_Santalo}
    E_\nu \zeta=\frac {\pi
    \abs{\mathcal{D}_1}} {\abs{v_1}\abs{\partial\mathcal{D}_1}}.
\end {equation}
If $\hat\zeta:\hat\Omega\rightarrow\mathbb{R} $ is the free
flight time between collisions with the piston, then it follows from
Proposition \ref{prop:inducing} that
\begin {equation}
\label{eq:2d_flight}
    E_{\hat\nu} \hat\zeta=\frac {\pi
    \abs{\mathcal{D}_1}}{\abs{v_1}\ell}.
\end{equation}

The expected value of $ \abs{v_1^\perp }$ when the left gas particle
collides with the (immobile) piston is given by
\begin {equation}
\label{eq:2d_momentum}
    E_{\hat\nu} \abs{v_1^\perp }=E_{\hat\nu} \sqrt{2E_1}\cos\varphi=
    \frac{\sqrt{2E_1}}{2}\int_{-\pi/2}^{+\pi/2} \cos^2\varphi\,d\varphi=
    \sqrt{2E_1}\frac{\pi}{4}.
\end {equation}

We wish to compute $\lim_{t\rightarrow\infty} t^{-1} \int_0^t
\abs{2v_1^\perp (s)}\delta_{q_1^\perp (s)=Q}ds$, the time average of
the change in momentum of the left gas particle when it collides
with the piston.  If this limit exists and is equal for almost every
initial condition of the left gas particle, then it makes sense to
define the pressure inside $\mathcal{D}_1$ to be this quantity
divided by $\ell$.  Because the collisions are hard-core, we cannot
directly apply Birkhoff's Ergodic Theorem to compute this limit.
However, we can compute this limit by using the map $\hat F$.

\begin {lem}
\label{lem:ae_convergence}

If the billiard flow $y(t) $ is ergodic, then for $\mu-a.e.$ $y\in
\mathcal{M}^1$,
\[
    \lim_{t\rightarrow\infty} \frac{1}{t}
    \int_0^t \abs{v_1^\perp (s)}\delta_{q_1^\perp (s)
    =Q}ds=
    \frac{E_1\ell}{2\abs{\mathcal{D}_1(Q)}}.
\]

\end {lem}

\begin {proof}
Because the billiard flow may be viewed as a suspension flow over
the collision cross-section with $\zeta$ as the height function, it
suffices to show that the convergence takes place for $\hat\nu-a.e.$
$x\in\hat\Omega$. For an initial condition $x\in\hat\Omega$, define
$\hat{N}_t(x)=\hat{N}_t=\#\set{s\in (0,t]:y(s,x) \in\hat\Omega}$. By
the Poincar\'e  Recurrence Theorem, $\hat{N}_t\rightarrow\infty$ as
$t\rightarrow\infty$, $\hat\nu-a.e.$

But
\[
\begin {split}
    \frac{\hat{N}_t}{\sum_{n=0}^{\hat{N}_t}\hat\zeta(\hat F^n x)}
    \frac{1}{\hat{N}_t}
    \sum_{n=1}^{\hat{N}_t}\abs{v_1^\perp }(\hat F^n x)
    &\leq
    \frac{1}{t}\int_0^t \abs{v_1^\perp (s)}\delta_{q_1^\perp (s)=Q}ds
    \\
    &\leq
    \frac{\hat{N}_t}{\sum_{n=0}^{\hat{N}_t-1}\hat\zeta(\hat F^n x)}
    \frac{1}{\hat{N}_t}
    \sum_{n=0}^{\hat{N}_t}\abs{v_1^\perp }(\hat F^n x),
\end {split}
\]
and so the result follows from Birkhoff's Ergodic Theorem and
Equations \eqref{eq:2d_flight} and \eqref{eq:2d_momentum}.
\end {proof}

\begin {cor}
\label{cor:ae_convergence}

If the billiard flow $y(t) $ is ergodic, then for each $\delta>0$,
\[
    \mu
    \set{y\in \mathcal{M}^1:\abs{\frac{1}{t}
    \int_0^t \abs{v_1^\perp (s)}\delta_{q_1^\perp (s)
    =Q}ds-
    \frac{E_1\ell}{2\abs{\mathcal{D}_1(Q)}}}\geq \delta}
    \rightarrow 0\text{ as }t\rightarrow \infty.
\]

\end {cor}

\subsection{Analysis of collisions}\label{sct:collisions}

In this section, we return to studying our piston system when
$\varepsilon>0$.  We will examine what happens when a particle
collides with the piston.  For convenience, we will only examine in
detail collisions between the piston and the left gas particle.
Collisions with the right gas particle can be handled similarly.

When the left gas particle collides with the piston, $v_1^\perp $
and $V$ instantaneously change according to the laws of elastic
collisions:
\begin {equation*}
    \begin{bmatrix}
    v_1^{\perp +}\\ V^+
    \end{bmatrix}
    =
    \frac{1}{1+M}
    \begin{bmatrix}
    1-M& 2M\\
    2& M-1\\
    \end{bmatrix}
    \begin{bmatrix}
    v_1^{\perp -}\\ V^-
    \end{bmatrix}.
\end {equation*}
In our coordinates, this becomes
\begin {equation}\label{eq:collision_change}
    \begin{bmatrix}
    v_1^{\perp +}\\ W^+
    \end{bmatrix}
    =
    \frac{1}{1+\varepsilon^2 }
    \begin{bmatrix}
    \varepsilon^2 -1 & 2\varepsilon\\
    2\varepsilon & 1-\varepsilon^2 \\
    \end{bmatrix}
    \begin{bmatrix}
    v_1^{\perp -}\\ W^-
    \end{bmatrix}.
\end {equation}
Recalling that $ v_1, W=\mathcal{O} (1) $, we find that to first
order in $\varepsilon$,
\begin{equation}
\label{eq:v_1Wchange}
    v_1^{\perp +}=-v_1^{\perp -}+\mathcal{O}(\varepsilon),\qquad
    W^ +=W^ -+\mathcal{O}(\varepsilon).
\end{equation}
Observe that a collision can only take place if $v_1^{\perp
-}>\varepsilon W^ - $.  In particular, $v_1^{\perp -}> -
\varepsilon\sqrt{2E_{max}}$.  Thus, either $v_1^{\perp -}> 0$ or
$v_1^{\perp -}= \mathcal{O} (\varepsilon) $.  By expanding
Equation~\eqref{eq:collision_change} to second order in
$\varepsilon$, it follows that
\begin{equation}
\label{eq:E_1Wchange}
\begin {split}
    E_1^+ -E_1^- &=-2\varepsilon W \abs{v_1^{\perp }}
    +\mathcal{O}(\varepsilon^2),\\
    W^+ -W^- &=+2\varepsilon \abs{v_1^{\perp }}
    +\mathcal{O}(\varepsilon^2).
\end {split}
\end{equation}
Note that it is immaterial whether we use the pre-collision or
post-collision values of $W$ and $\abs{v_1^{\perp }}$ on the right
hand side of Equation~\eqref{eq:E_1Wchange}, because any ambiguity
can be absorbed into the $\mathcal{O} (\varepsilon^2) $ term.

It is convenient for us to define a ``clean collision'' between the
piston and the left gas particle:
\begin {defn}
    The left gas particle experiences a \emph{clean collision} with the
    piston if and only if $v_1^{\perp -}>0$ and $v_1^{\perp +}<-\varepsilon
    \sqrt{2E_{max}}$.
\end{defn}
\noindent In particular, after a clean collision, the left gas
particle will escape from the piston, i.e.~the left gas particle
will have to move into the region $q_1^{\perp }\leq 0 $ before it
can experience another collision with the piston.  It follows that
there exists a constant $C_1>0$, which depends on the set
$\mathcal{V} $, such that for all $\varepsilon$ sufficiently small,
so long as $ Q\geq Q_{min} $ and $\abs{v_1^{\perp }}
>\varepsilon C_1$ when $q_1^{\perp }\in [Q_{ min},Q]$, then the left gas particle
will experience only clean collisions with the piston, and the time
between these collisions will be greater than $2Q_{min}/(\sqrt
{2E_{max}})$.  (Note that when we write expressions such as
$q_1^{\perp }\in [Q_{ min},Q]$, we implicitly mean that $q_1$ is
positioned inside the ``tube'' discussed at the beginning of
Section~\ref{sct:main_result}.) One can verify that $C_1=5\sqrt
{2E_{max}}$ would work.

Similarly, we can define clean collisions between the right gas
particle and the piston.  We assume that $C_1$ was chosen
sufficiently large such that for all $\varepsilon$ sufficiently
small, so long as $ Q\leq Q_{max} $ and $\abs{v_2^{\perp }}
>\varepsilon C_1$ when $q_2^{\perp }\in [Q,Q_{max}] $, then the right gas particle
will experience only clean collisions with the piston.

Now we define three more stopping times, which are functions of the
initial conditions in phase space.
\[
\begin {split}
    T_\varepsilon' =&\inf \{\tau\geq 0: Q_{min}\leq q_{1,\varepsilon}^{\perp
    }(\tau/\varepsilon)\leq Q_\varepsilon(\tau/\varepsilon)\leq Q_{max}
    \text { and}\abs{v_{1,\varepsilon}^{\perp }(\tau/\varepsilon)}\leq C_1\varepsilon \}
    ,
    \\
    T_\varepsilon'' =&
    \inf \{\tau\geq 0:Q_{min}\leq Q_\varepsilon(\tau/\varepsilon)\leq q_{2,\varepsilon}^{\perp
    }(\tau/\varepsilon)\leq Q_{max}
    \text { and}\abs{v_{2,\varepsilon}^{\perp }(\tau/\varepsilon)}\leq C_1\varepsilon
    \},
    \\
    \tilde{T}_\varepsilon =&
    T\wedge T_\varepsilon\wedge T_\varepsilon'\wedge T_\varepsilon''
\end {split}
\]

Define $H(z) $ by
\[
    H(z) =
    \begin{bmatrix}
    W\\
    +2\abs{v_1^{\perp }} \delta_{q_1^{\perp }=Q}
    -2\abs{v_2^{\perp }} \delta_{q_2^{\perp }=Q}\\
    -2W\abs{v_1^{\perp }} \delta_{q_1^{\perp }=Q}\\
    +2W\abs{v_2^{\perp }} \delta_{q_2^{\perp }=Q}\\
    \end{bmatrix}.
\]
Here we make use of Dirac delta functions.  All integrals involving
these delta functions may be replaced by sums.

The following lemma is an immediate consequence of Equation
\eqref{eq:E_1Wchange} and the above discussion:

\begin{lem}
\label{lem:h_int}

If $0\leq t_1\leq t_2\leq \tilde{T}_\varepsilon/\varepsilon $, the
piston experiences $\mathcal{O} ((t_2-t_1)\vee 1) $ collisions with
gas particles in the time interval $[t_1, t_2]$, all of which are
clean collisions. Furthermore,
\begin {equation*}
    h_\varepsilon(t_2)-h_\varepsilon(t_1)=
    \mathcal{O}(\varepsilon)+\varepsilon\int_{t_1}^{t_2}
    H(z_\varepsilon(s))ds.
\end {equation*}
Here any ambiguities arising from collisions occurring at the limits
of integration can be absorbed into the $\mathcal{O} (\varepsilon) $
term.

\end {lem}

\subsection{Another heuristic derivation of the averaged
equation}\label{sct:heuristic2}

The following heuristic derivation of Equation \eqref{eq:d_davg}
when $ d=2$ was suggested in~\cite{Dol05}. Let $\Delta t $ be a
length of time long enough such that the piston experiences many
collisions with the gas particles, but short enough such that the
slow variables change very little, in this time interval.  From each
collision with the left gas particle, Equation~\eqref{eq:E_1Wchange}
states that $W$ changes by an amount $+2\varepsilon
\abs{v_1^{\perp}} +\mathcal{O}(\varepsilon^2)$, and from
Equation~\eqref{eq:2d_momentum} the average change in $W$ at these
collisions should be approximately $\varepsilon\pi
\sqrt{2E_1}/2+\mathcal{O}(\varepsilon^2)$. From
Equation~\eqref{eq:2d_flight} the frequency of these collisions is
approximately $\sqrt{2E_1}\,\ell /(\pi
    \abs{\mathcal{D}_1})$. Arguing
similarly for collisions with the other particle, we guess that
\[
    \frac {\Delta W} {\Delta t} =
    \varepsilon\frac{E_1\ell}{\abs{\mathcal{D}_1(Q)}}
    -\varepsilon\frac{E_2\ell}{\abs{\mathcal{D}_2(Q)}}
    +\mathcal{O}(\varepsilon^2).
\]
With $\tau=\varepsilon t$ as the slow time, a reasonable guess for
the averaged equation for $W$ is
\[
    \frac {dW}{d\tau}=\frac{E_1\ell}{\abs{\mathcal{D}_1(Q)}}
    -\frac{E_2\ell}{\abs{\mathcal{D}_2(Q)}}.
\]
Similar arguments for the other slow variables lead to the averaged
equation \eqref{eq:d_davg}, and this explains why we used $P_i=
E_i/\abs{\mathcal{D}_i}$ for the pressure of a $2$-dimensional gas
in Section~\ref{sct:heuristic}.

There is a similar heuristic derivation of the averaged equation in
$ d>2$ dimensions.  Compare the analogues of
Equations~\eqref{eq:2d_flight} and \eqref{eq:2d_momentum} in
Subsection~\ref{sct:higher_d}.

\subsection{\textit{A priori} estimate on the size
        of a set of bad initial conditions}

In this section, we give an \textit{a priori} estimate on the size
of a set of initial conditions that should not give rise to orbits
for which $\sup_{0\leq\tau\leq T\wedge
T_\varepsilon}\abs{h_\varepsilon(\tau/\varepsilon)-\bar{h}(\tau)}$
is small.  In particular, when proving Theorem \ref{thm:dDpiston},
it is convenient to focus on orbits that only contain clean
collisions with the piston.  Thus, we show that $P\{\tilde
{T}_\varepsilon<T\wedge T_\varepsilon \} $ vanishes as
$\varepsilon\rightarrow 0$. At first, this result may seem
surprising, since $P\{T_\varepsilon'\wedge T_\varepsilon''=0\}
=\mathcal{O}(\varepsilon)$, and one would expect $\cup_{t=0} ^
{T/\varepsilon} z_\varepsilon(-t)\{T_\varepsilon'\wedge
T_\varepsilon''=0\}$ to have a size of order $1$. However, the rate
at which orbits escape from $\{T_\varepsilon'\wedge
T_\varepsilon''=0\}$ is very small, and so we can prove the
following:

\begin {lem}
\label{lem:no_vertical}
    \[
        P\{\tilde {T}_\varepsilon<T\wedge T_\varepsilon \} =\mathcal{O}
        (\varepsilon).
    \]
\end {lem}

In some sense, this lemma states that the probability of having a
gas particle move nearly parallel to the piston's face within the
time interval $[0,T/\varepsilon ] $, when one would expect the other
gas particle to force the piston to move on a macroscopic scale,
vanishes as $\varepsilon\rightarrow 0$.  Thus, one can hope to
control the occurrence of the ``nondiffusive fluctuations'' of the
piston described in~\cite{CD06} on a time scale $\mathcal{O}
(\varepsilon^ {-1}) $.

\begin {proof}

As the left and the right gas particles can be handled similarly, it
suffices to show that $P\{T_\varepsilon'<T \} =\mathcal{O}
(\varepsilon)$.  Define
\[
    \mathfrak{B}_\varepsilon=\{z\in\mathcal{M}:
    Q_{min}\leq q_1^{\perp
    }\leq Q\leq Q_{max}
    \text { and}\abs{v_1^{\perp }}\leq C_1\varepsilon
    \}.
\]
Then $\{T_\varepsilon'<T \}\subset \cup_{t=0} ^ {T/\varepsilon}
z_\varepsilon(-t)\mathfrak{B}_\varepsilon$, and if $\gamma=
Q_{min}/\sqrt{8 E_{max}}$,
\[
\begin {split}
    P\left (\bigcup_{t=0} ^{T/\varepsilon}
    z_\varepsilon(-t)\mathfrak{B}_\varepsilon\right)
    &
    =
    P\left (\bigcup_{t=0} ^{T/\varepsilon}
    z_\varepsilon(t)\mathfrak{B}_\varepsilon\right)
    =
    P\left(\mathfrak{B}_\varepsilon\cup\bigcup_{t=0} ^{T/\varepsilon}
    ((z_\varepsilon(t)\mathfrak{B}_\varepsilon) \backslash \mathfrak{B}_\varepsilon)
    \right)
    \\
    &
    \leq
    P\mathfrak{B}_\varepsilon+P\left( \bigcup_{k=0}^{T/(\varepsilon\gamma
    )} z_\varepsilon(k\gamma )
    \Bigl[ \bigcup_{t=0}^\gamma
    (z_\varepsilon(t)\mathfrak{B}_\varepsilon)\backslash
    \mathfrak{B}_\varepsilon \Bigr]
    \right)
    \\
    &
    \leq
    P\mathfrak{B}_\varepsilon+
    \left(\frac{T}{\varepsilon\gamma }+1\right)
    P\left(
    \bigcup_{t=0}^\gamma
    (z_\varepsilon(t)\mathfrak{B}_\varepsilon)\backslash
    \mathfrak{B}_\varepsilon
    \right).
\end {split}
\]
Now $P\mathfrak{B}_\varepsilon=\mathcal{O} (\varepsilon) $, so if we
can show that
$P\left(\bigcup_{t=0}^\gamma(z_\varepsilon(t)\mathfrak{B}_\varepsilon)\backslash
\mathfrak{B}_\varepsilon\right)=\mathcal{O} (\varepsilon^2)$, then
it will follow that $P\{T_\varepsilon'<T \} =\mathcal{O}
(\varepsilon)$.

If
$z\in\bigcup_{t=0}^\gamma(z_\varepsilon(t)\mathfrak{B}_\varepsilon)\backslash
\mathfrak{B}_\varepsilon$, it is still true that
$\abs{v_1^\perp}=\mathcal{O}(\varepsilon)$.  This is because
$\abs{v_1^\perp}$ changes by at most $\mathcal{O} (\varepsilon) $ at
the collisions, and if a collision forces
$\abs{v_1^\perp}>C_1\varepsilon$, then the gas particle must escape
to the region $q_1^\perp\leq 0$ before $ v_1^\perp $ can change
again, and this will take time greater than $\gamma $. Furthermore,
if
$z\in\bigcup_{t=0}^\gamma(z_\varepsilon(t)\mathfrak{B}_\varepsilon)\backslash
\mathfrak{B}_\varepsilon$, then at least one of the following four
possibilities must hold:
\begin {itemize}
\item
     $\abs{q_1^\perp-Q_{min}}\leq\mathcal{O}(\varepsilon)$,
\item
     $\abs{Q-Q_{min}}\leq\mathcal{O}(\varepsilon)$,
\item
     $\abs{Q-Q_{max}}\leq\mathcal{O}(\varepsilon)$,
\item
     $\abs{Q-q_1^\perp}\leq\mathcal{O}(\varepsilon)$.
\end {itemize}
It follows that
$P\left(\bigcup_{t=0}^\gamma(z_\varepsilon(t)\mathfrak{B}_\varepsilon)\backslash
\mathfrak{B}_\varepsilon\right)=\mathcal{O} (\varepsilon^2)$.  For
example,
\[
\begin {split}
    \int_{\mathcal{M}}
    &
    1_{\{\abs{v_1^\perp}\leq
    \mathcal{O}(\varepsilon),\:
    \abs{q_1^\perp-Q_{min}}\leq\mathcal{O}(\varepsilon)\}}dP
    \\
    &
    =
    \C
    \int_{\set {E_{min}\leq W^2/2+v_1^2/2+v_2^2/2\leq E_{max}}}
    1_{\{\abs{v_1^\perp}\leq
    \mathcal{O}(\varepsilon)\}}
    dW dv_1^{\perp}dv_1^{\parallel} dv_2^{\perp}dv_2^{\parallel}
    \\
    &
    \qquad\times
    \int_{\set {Q\in [0,1],\, q_1\in\mathcal{D}_1,\, q_2\in\mathcal{D}_2}}
    1_{\{\abs{q_1^\perp-Q_{min}}\leq\mathcal{O}(\varepsilon)\}}
    dQ dq_1^{\perp}dq_1^{\parallel} dq_2^{\perp}dq_2^{\parallel}
    \\
    &
    =\mathcal{O}(\varepsilon^2).
\end {split}
\]

\end {proof}

\section{Proof of the main result for two-dimensional gas
    containers with only one gas particle on each side}
\label{sct:2dproof}

As in Section~\ref{sct:2dprep}, we continue with the case when $d=2$
and there is only one gas particle on either side of the piston.

\subsection{Main steps in the proof of convergence in probability}\label{sct:main_steps}

By Lemma \ref{lem:no_vertical}, it suffices to show that $
\sup_{0\leq\tau\leq \tilde{T}_\varepsilon}
\abs{h_\varepsilon(\tau/\varepsilon)-\bar{h}(\tau)}\rightarrow 0$ in
probability as $\varepsilon=M^{-1/2}\rightarrow 0$. Several of the
ideas in the steps below were inspired by a recent proof of Anosov's
averaging theorem for smooth systems that is due to
Dolgopyat~\cite{Dol05}.

\paragraph*{Step 1:  Reduction using Gronwall's Inequality.}

Observe that $\bar {h}(\tau) $ satisfies the integral equation
\[
\bar {h}(\tau) -\bar h(0) = \int_0^{\tau}\bar H(\bar
h(\sigma))d\sigma,
\]
while from Lemma \ref{lem:h_int},
\[
\begin {split}
    h_\varepsilon(\tau/\varepsilon)-h_\varepsilon(0)
    &
    =\mathcal{O}(\varepsilon) +\varepsilon\int_0^{\tau/\varepsilon}
    H(z_\varepsilon(s))ds\\
    &=\mathcal{O}(\varepsilon) +
    \varepsilon\int_0^{\tau/\varepsilon}
    H(z_\varepsilon(s))-
    \bar H(h_\varepsilon(s))ds+
    \int_0^{\tau}\bar H( h_\varepsilon(\sigma/\varepsilon))d\sigma
\end {split}
\]
for $0\leq\tau\leq \tilde{T}_\varepsilon$.  Define
\[
    e_\varepsilon(\tau) =\varepsilon\int_0^{\tau/\varepsilon}
    H(z_\varepsilon(s))- \bar H(h_\varepsilon(s))ds.
\]
It follows from
Gronwall's Inequality that
\begin {equation}\label {eq:2d_Gronwall}
    \sup_{0\leq \tau\leq \tilde{T}_\varepsilon}
    \abs{h_\varepsilon(\tau/\varepsilon)-\bar h(\tau)}\leq
    \left(\mathcal{O}(\varepsilon)+
    \sup_{0\leq \tau\leq \tilde{T}_\varepsilon}
    \abs{e_\varepsilon(\tau)}\right)e^{ \Lip{\bar
    H\arrowvert_\mathcal{V}}T}.
\end {equation}
\noindent Gronwall's Inequality is usually stated for continuous
paths, but the standard proof (found in \cite{SV85}) still works for
paths that are merely integrable, and
$\abs{h_\varepsilon(\tau/\varepsilon)-\bar h(\tau)}$ is piecewise
smooth.

\paragraph*{Step 2:  Introduction of a time scale for ergodization.}

Let $L(\varepsilon) $ be a real valued function such that
$L(\varepsilon)\rightarrow\infty$, but $L(\varepsilon)\ll \log
\varepsilon^ {-1} $, as $\varepsilon\rightarrow 0$. In
Section~\ref{sct:Gronwall} we will place precise restrictions on the
growth rate of $L(\varepsilon) $.  Think of $L(\varepsilon) $ as
being a time scale that grows as $\varepsilon\rightarrow 0$ so that
\emph{ergodization}, i.e.~the convergence along an orbit of a
function's time average to a space average, can take place. However,
$L(\varepsilon) $ doesn't grow too fast, so that on this time scale
$z_\varepsilon(t) $ essentially stays on the submanifold
$\set{h=h_\varepsilon(0)}$, where we have our ergodicity assumption.
Set $t_{k,\varepsilon} =kL(\varepsilon) $, so that
\begin {equation}
\label{eq:2d_einfnorm}
    \sup_{0\leq \tau\leq \tilde{T}_\varepsilon}\abs{e_\varepsilon(\tau)}
    \leq \mathcal{O}(\varepsilon L(\varepsilon))+
    \varepsilon\sum_{k=0}^{\frac{\tilde{T}_\varepsilon}{\varepsilon
    L(\varepsilon)}-1}\abs{\int_{t_{k,\varepsilon}}^{t_{k+1,\varepsilon}}H(z_\varepsilon(s))-\bar
    H(h_\varepsilon(s))ds}.
\end {equation}

\paragraph*{Step 3:  A splitting according to particles.}

Now $H(z) -\bar H(h(z)) $ divides into two pieces, each of which
depends on only one gas particle when the piston is held fixed:
\[
    H(z) -\bar H(h(z))=
    \begin{bmatrix}
    0\\
    2\abs{v_1^{\perp }} \delta_{q_1^{\perp }=Q}
    -\frac{E_1\ell}{\abs{\mathcal{D}_1(Q)}}\\
    -2W\abs{v_1^{\perp }} \delta_{q_1^{\perp }=Q}
    +\frac{WE_1\ell}{\abs{\mathcal{D}_1(Q)}}\\
    0\\
    \end{bmatrix}
    +
    \begin{bmatrix}
    0\\
    \frac{E_2\ell}{\abs{\mathcal{D}_2(Q)}}
    -2\abs{v_2^{\perp }} \delta_{q_2^{\perp }=Q}\\
    0\\
    -\frac{WE_2\ell}{\abs{\mathcal{D}_2(Q)}}+
    2W\abs{v_2^{\perp }} \delta_{q_2^{\perp }=Q}\\
    \end{bmatrix}.
\]
We will only deal with the piece depending on the left gas particle,
as the right particle can be handled similarly.  Define
\begin{equation}
\label{eq:G_definition}
    G(z)=\abs{v_1^{\perp }} \delta_{q_1^{\perp }=Q},
    \qquad
    \bar G(h)=
    \frac{E_1\ell}{2\abs{\mathcal{D}_1(Q)}}.
\end {equation}
Returning to Equation \eqref{eq:2d_einfnorm}, we see that in order
to prove Theorem \ref{thm:dDpiston}, it suffices to show that both
\[
\begin {split}
    &\varepsilon\sum_{k=0}^{\frac{\tilde{T}_\varepsilon}{\varepsilon
    L(\varepsilon)}-1}\abs{\int_{t_{k,\varepsilon}}^{t_{k+1,\varepsilon}}G(z_\varepsilon(s))-\bar
    G(h_\varepsilon(s))ds} \text { and}
    \\
    &\varepsilon\sum_{k=0}^{\frac{\tilde{T}_\varepsilon}{\varepsilon
    L(\varepsilon)}-1}\abs{\int_{t_{k,\varepsilon}}^{t_{k+1,\varepsilon}}W_\varepsilon(s)
    \bigl(G(z_\varepsilon(s))-\bar
    G(h_\varepsilon(s))\bigr)ds}
\end {split}
\]
converge to $0$ in probability as $\varepsilon\rightarrow 0$.

\paragraph*{Step 4:  A splitting for using the triangle inequality.}

Now we let $z_{k,\varepsilon} (s) $ be the orbit of the
$\varepsilon=0$ Hamiltonian vector field satisfying
$z_{k,\varepsilon}(t_{k,\varepsilon})=z_{\varepsilon}(t_{k,\varepsilon})$.
Set $h_{k,\varepsilon} (t) =h (z_{k,\varepsilon} (t)) $. Observe
that $h_{k,\varepsilon} (t) $ is independent of $t $.

We emphasize that so long as $0\leq
t\leq\tilde{T}_\varepsilon/\varepsilon$, the times between
collisions of a specific gas particle and piston are uniformly
bounded greater than $0$, as explained before Lemma \ref{lem:h_int}.
It follows that, so long as
$t_{k+1,\varepsilon}\leq\tilde{T}_\varepsilon/\varepsilon$,
\begin{equation}
\label{eq:h_div}
    \sup_{t_{k,\varepsilon}\leq t\leq t_{k+1,\varepsilon}}
    \abs{h_{k,\varepsilon} (t) -h_\varepsilon (t)}
    =\mathcal{O}(\varepsilon L(\varepsilon)).
\end{equation}
This is because the slow variables change by at most $\mathcal{O}
(\varepsilon) $ at collisions, and
$dQ_\varepsilon/dt=\mathcal{O}(\varepsilon)$.

Also,
\[
\begin {split}
    \int_{t_{k,\varepsilon}}^{t_{k+1,\varepsilon}}
    &
    W_\varepsilon(s) \bigl(G(z_\varepsilon(s))-\bar
    G(h_\varepsilon(s))\bigr)ds
    \\
    &= \mathcal{O} (\varepsilon
    L(\varepsilon)^2) +W_{k,\varepsilon}(t_{k,\varepsilon})\int_{t_{k,\varepsilon}}^{t_{k+1,\varepsilon}}
    G(z_\varepsilon(s))-\bar G(h_\varepsilon(s))ds,
\end {split}
\]
and so
\[
\begin {split}
    \varepsilon\sum_{k=0}^{\frac{\tilde{T}_\varepsilon}{\varepsilon
    L(\varepsilon)}-1}
    &
    \abs{\int_{t_{k,\varepsilon}}^{t_{k+1,\varepsilon}}W_\varepsilon(s)
    \bigl(G(z_\varepsilon(s))-\bar
    G(h_\varepsilon(s))\bigr)ds}
    \\
    &\leq
    \mathcal{O}(\varepsilon L(\varepsilon))+
    \varepsilon\,\C\sum_{k=0}^{\frac{\tilde{T}_\varepsilon}{\varepsilon
    L(\varepsilon)}-1}\abs{\int_{t_{k,\varepsilon}}^{t_{k+1,\varepsilon}}
    G(z_\varepsilon(s))-\bar
    G(h_\varepsilon(s))ds}.
\end {split}
\]
Thus, in order to prove Theorem \ref{thm:dDpiston}, it suffices to
show that
\[
\begin {split}
    \varepsilon\sum_{k=0}^{\frac{\tilde{T}_\varepsilon}{\varepsilon
    L(\varepsilon)}-1}\abs{\int_{t_{k,\varepsilon}}^{t_{k+1,\varepsilon}}G(z_\varepsilon(s))-\bar
    G(h_\varepsilon(s))ds}
    \leq
    \varepsilon\sum_{k=0}^{\frac{\tilde{T}_\varepsilon}{\varepsilon
    L(\varepsilon)}-1}
    \abs{I_{k,\varepsilon}}+\abs{II_{k,\varepsilon}}+\abs{III_{k,\varepsilon}}
\end {split}
\]
converges to $0$ in probability as $\varepsilon\rightarrow 0$, where
\[
\begin {split}
    I_{k,\varepsilon}
    &=\int_{t_{k,\varepsilon}}^{t_{k+1,\varepsilon}}G(z_\varepsilon(s)) -
    G(z_{k,\varepsilon}(s))ds,
    \\
    II_{k,\varepsilon}
    &=\int_{t_{k,\varepsilon}}^{t_{k+1,\varepsilon}} G(z_{k,\varepsilon}(s))-
    \bar G(h_{k,\varepsilon}(s))ds,
    \\
    III_{k,\varepsilon}
    &=\int_{t_{k,\varepsilon}}^{t_{k+1,\varepsilon}} \bar G(h_{k,\varepsilon}(s))-
    \bar G(h_{\varepsilon}(s))ds.
\end {split}
\]

The term $II_{k,\varepsilon}$ represents an ``ergodicity term'' that
can be controlled by our assumptions on the ergodicity of the flow
$z_0(t) $, while the terms $I_{k,\varepsilon}$ and
$III_{k,\varepsilon}$ represent ``continuity terms'' that can be
controlled by controlling the drift of $z_{\varepsilon} (t) $ from
$z_{k,\varepsilon} (t) $ for $t_{k,\varepsilon}\leq t\leq
t_{k+1,\varepsilon}$.

\paragraph*{Step 5:  Control of drift from the $\varepsilon=0$ orbits.}

Now $\bar G$ is uniformly Lipschitz on the compact set $\mathcal{V}
$, and so it follows from Equation \eqref{eq:h_div} that
$III_{k,\varepsilon}=\mathcal{O}(\varepsilon L(\varepsilon)^2)$.
Thus,
$\varepsilon\sum_{k=0}^{\frac{\tilde{T}_\varepsilon}{\varepsilon
L(\varepsilon)}-1}\abs{III_{k,\varepsilon}} =\mathcal{O}(\varepsilon
L(\varepsilon))\rightarrow 0$ as $\varepsilon\rightarrow 0$.

Next, we show that for fixed $\delta > 0$,
$P\left(\varepsilon\sum_{k=0}^{\frac{\tilde{T}_\varepsilon}{\varepsilon
    L(\varepsilon)}-1}\abs{I_{k,\varepsilon}}\geq
\delta\right)\rightarrow 0$ as $\varepsilon\rightarrow 0$.

For initial conditions $z\in \mathcal{M}$ and for integers $k\in
[0,T/(\varepsilon L(\varepsilon))-1]$ define
\[
\begin {split}
    \mathcal{A}_{k,\varepsilon} & =\set{z:\frac{1}{L(\varepsilon)}
    \abs{I_{k,\varepsilon}}
    >\frac{\delta}{2T} \text { and } k\leq\frac{\tilde{T}_\varepsilon}{\varepsilon
    L(\varepsilon)}-1 }
    ,\\
    \mathcal{A}_{z,\varepsilon} & =\set{k:z\in
    \mathcal{A}_{k,\varepsilon}}.
\end {split}
\]
Think of these sets as describing ``poor continuity'' between
solutions of the $\varepsilon=0$ and the $\varepsilon>0 $
Hamiltonian vector fields. For example, roughly speaking,
$z\in\mathcal{A}_{k,\varepsilon}$ if the orbit $z_{\varepsilon}(t)$
starting at $z$ does not closely follow $z_{k,\varepsilon}(t)$ for
$t_{k,\varepsilon}\leq t\leq t_{k+1,\varepsilon} $.

One can easily check that $\abs{I_{k,\varepsilon}}\leq \mathcal{O}
(L(\varepsilon))$ for $k\leq\ \tilde{T}_\varepsilon/(\varepsilon
L(\varepsilon))-1$, and so it follows that
\[
    \varepsilon
    \sum_{k=0}^{\frac{\tilde{T}_\varepsilon}{\varepsilon
    L(\varepsilon)}-1}\abs{I_{k,\varepsilon}}\leq
    \frac{\delta}{2}+\mathcal{O}(\varepsilon L(\varepsilon)
    \# (\mathcal{A}_{z,\varepsilon})).
\]
Therefore it suffices to show that $P\left(\#
(\mathcal{A}_{z,\varepsilon})\geq\delta(\C\, \varepsilon
L(\varepsilon)) ^ {-1}\right)\rightarrow 0$ as
$\varepsilon\rightarrow 0$. By Chebyshev's Inequality, we need only
show that
\[
    E_P (\varepsilon L(\varepsilon)\# (\mathcal{A}_{z,\varepsilon})) =
    \varepsilon L(\varepsilon)\sum_{k=0}^{\frac{T}{\varepsilon
    L(\varepsilon)}-1}P(\mathcal{A}_{k,\varepsilon})
\]
tends to $0$ with $\varepsilon$.

Observe that $z_\varepsilon
(t_{k,\varepsilon})\mathcal{A}_{k,\varepsilon}\subset\mathcal{A}_{0,\varepsilon}
$. In words, the initial conditions giving rise to orbits that are
``bad'' on the time interval
$[t_{k,\varepsilon},t_{k+1,\varepsilon}] $, moved forward by time
$t_{k,\varepsilon}$, are initial conditions giving rise to orbits
which are ``bad'' on the time interval
$[t_{0,\varepsilon},t_{1,\varepsilon}] $. Because the flow
$z_\varepsilon(\cdot) $ preserves the measure, we find that
\[
    \varepsilon L(\varepsilon)\sum_{k=0}^{\frac{T}{\varepsilon
    L(\varepsilon)}-1}P(\mathcal{A}_{k,\varepsilon})
    \leq \C\, P(\mathcal{A}_{0,\varepsilon}).
\]

To estimate $P(\mathcal{A}_{0,\varepsilon})$, it is convenient to
use a different probability measure, which is uniformly equivalent
to $P$ on the set $\set{z\in\mathcal{M}:h(z)\in\mathcal{V}}
\supset\{\tilde{T}_\varepsilon\geq \varepsilon L(\varepsilon)\} $.
We denote this new probability measure by $P^f$, where the $f$
stands for ``factor.''  If we choose coordinates on $\mathcal{M} $
by using $h$ and the billiard coordinates on the two gas particles,
then $P^f$ is defined on $\mathcal{M}$ by $dP^f
=dh\,d\mu^1_h\,d\mu^2_h$, where $dh$ represents the uniform measure
on $\mathcal{V}\subset\mathbb{R}^4$, and the factor measure
$d\mu^i_h$ represents the invariant billiard measure of the $i^ {th}
$ gas particle coordinates for a fixed value of the slow variables.
One can verify that $1_{\set{h(z)\in\mathcal{V}}}dP \leq \C\, dP^f$,
but that $P^f$ is not invariant under the flow $z_\varepsilon(\cdot)
$ when $\varepsilon>0$.

We abuse notation, and consider $\mu^1_h$ to be a measure on the
left particle's initial billiard coordinates once $h$ and the
initial coordinates of the right gas particle are fixed.  In this
context, $\mu^1_h$ is simply the measure $\mu$ from
Subsection~\ref{sct:billiard}.  Then
\[
\begin {split}
    &
    P^f(\mathcal{A}_{0,\varepsilon})
    \\
    &\leq \int
    dh\,d\mu^2_h \cdot\mu_h^1
    \set{z:\abs{\frac{1}{L(\varepsilon)}\int_0^{L(\varepsilon)}
    G(z_\varepsilon(s))-G(z_0(s))ds}\geq\frac{\delta}{2T}
    \text { and } \varepsilon
    L(\varepsilon)\leq\tilde{T}_\varepsilon },
\end {split}
\]
and we must show that the last term tends to $0$ with $\varepsilon$.
By the Bounded Convergence Theorem, it suffices to show that for
almost every $h\in\mathcal{V}$ and initial condition for the right
gas particle,
\begin {equation}\label{eq:Gronwall_probability}
    \mu_h^1
    \set{z:\abs{\frac{1}{L(\varepsilon)}\int_0^{L(\varepsilon)}
    G(z_\varepsilon(s))-G(z_0(s))ds}\geq\frac{\delta}{2T}
    \text { and } \varepsilon
    L(\varepsilon)\leq\tilde{T}_\varepsilon }
    \rightarrow 0\text{ as }\varepsilon\rightarrow 0.
\end {equation}

Note that if $G$ were a smooth function and $z_\varepsilon(\cdot) $
were the flow of a smooth family of vector fields $Z(z,\varepsilon)
$ that depended smoothly on $\varepsilon$, then from Gronwall's
Inequality, it would follow that $\sup_{0\leq t\leq
L(\varepsilon)}\abs{z_{\varepsilon}(t)-
    z_0(t)}\leq \mathcal{O}(\varepsilon L(\varepsilon)
    e^{ \Lip{Z} L(\varepsilon)}).$
If this were the case, then $\abs{L(\varepsilon)^
{-1}\int_0^{L(\varepsilon)}
G(z_\varepsilon(s))-G(z_0(s))ds}=\mathcal{O}(\varepsilon
L(\varepsilon) e^{ \Lip{Z} L(\varepsilon)})$, which would tend to
$0$ with $\varepsilon$.  Thus, we need a Gronwall-type inequality
for billiard flows.  We obtain the appropriate estimates in Section
\ref{sct:Gronwall}.

\paragraph*{Step 6:  Use of ergodicity along fibers to
control $II_{k,\varepsilon} $.}

All that remains to be shown is that for fixed $\delta > 0$,
$P\left(\varepsilon\sum_{k=0}^{\frac{\tilde{T}_\varepsilon}{\varepsilon
    L(\varepsilon)}-1}\abs{II_{k,\varepsilon}}\geq
\delta\right)\rightarrow 0$ as $\varepsilon\rightarrow 0$.

For initial conditions $z\in \mathcal{M}$ and for integers $k\in
[0,T/(\varepsilon L(\varepsilon))-1]$ define
\[
\begin {split}
    \mathcal{B}_{k,\varepsilon} & =\set{z:\frac{1}{L(\varepsilon)}
    \abs{II_{k,\varepsilon}}
    >\frac{\delta}{2T} \text { and } k\leq\frac{\tilde{T}_\varepsilon}{\varepsilon
    L(\varepsilon)}-1 }
    ,\\
    \mathcal{B}_{z,\varepsilon} & =\set{k:z\in
    \mathcal{B}_{k,\varepsilon}}.
\end {split}
\]
Think of these sets as describing ``bad ergodization.''  For
example, roughly speaking, $z\in\mathcal{B}_{k,\varepsilon}$ if the
orbit $z_{\varepsilon}(t)$ starting at $z$ spends the time between
$t_{k,\varepsilon}$ and $t_{k+1,\varepsilon} $ in a region of phase
space where the function $G(\cdot) $ is ``poorly ergodized'' on the
time scale $L(\varepsilon) $ by the flow $z_0(t) $ (as measured by
the parameter $\delta/2T$).  Note that $G(z)=\abs{v_1^{\perp }}
\delta_{q_1^{\perp }=Q}$ is not really a function, but that we may
still speak of the convergence of $t^ {-1}\int_0^t G(z_0(s))ds$ as
$t\rightarrow\infty$.  As we showed in
Lemma~\ref{lem:ae_convergence}, the limit is $\bar G(h_0) $ for
almost every initial condition.

Proceeding as in Step 5 above, we find that it suffices to show that
for almost every $h\in\mathcal{V}$,
\[
    \mu_h^1
    \set{z:\abs{\frac{1}{t}\int_0^{t}
    G(z_0(s))ds-\bar G(h_0(0))}\geq\frac{\delta}{2T}}
    \rightarrow 0\text{ as }t\rightarrow \infty.
\]
But this is simply a question of examining billiard flows, and it
follows immediately from Corollary \ref{cor:ae_convergence} and our
Main Assumption.

\subsection{A Gronwall-type inequality for billiards}\label{sct:Gronwall}

We begin by presenting a general version of Gronwall's Inequality
for billiard maps.  Then we will show how these results imply the
convergence required in Equation~\eqref{eq:Gronwall_probability}.

\subsubsection{Some inequalities for the collision map}\label{sct:Gronwall_map}

In this section, we consider the value of the slow variables to be
fixed at $h_0\in\mathcal{V} $.  We will use the notation and results
presented in Section~\ref{sct:billiard}, but because the value of
the slow variables is fixed, we will omit it in our notation.

Let $\rho $, $\gamma $, and $\lambda$ satisfy $0<\rho\ll\gamma \ll
1\ll \lambda<\infty $. Eventually, these quantities will be chosen
to depend explicitly on $\varepsilon$, but for now they are fixed.

Recall that the phase space $\Omega$ for the collision map $F$ is a
finite union of disjoint rectangles and cylinders.  Let $d$ be the
Euclidean metric on connected components of $\Omega$.  If $x$ and
$x'$ belong to different components, then we set $d(x,x') =\infty $.
The invariant measure $\nu$ satisfies $\nu<\C\cdot (\text {Lebesgue
measure}) $.  For $A\subset\Omega$ and $a >0$, let $\mathcal{N}_a
(A) =\set {x\in\Omega:d(x,A)< a} $ be the $a$-neighborhood of $A$.

For $x\in\Omega$ let $x_k(x) = x_k =F^k x$, $ k\geq 0$, be its
forward orbit.  Suppose $x\notin \mathcal{C}_{\gamma ,\lambda}$,
where
\[
    \mathcal{C}_{\gamma,\lambda}=
    \bigl(\cup_{k=0} ^\lambda F^
    {-k}\mathcal{N}_\gamma (\partial\Omega)\bigr)\bigcup\bigl(\cup_{k=0}
    ^\lambda F^ {-k}\mathcal{N}_\gamma  (F^ {-1}\mathcal{N}_\gamma
    (\partial\Omega))\bigr).
\]
Thus for $0\leq k\leq \lambda $, $x_k$ is well defined, and from
Equation~\eqref{eq:2d_derivative_bound} it satisfies
\begin {equation}\label{eq:gron0}
    d(x',x_k)\leq \gamma \;\Rightarrow\; d(Fx',x_{k+1})\leq\frac{\C}{\gamma
    } d(x',x_k).
\end {equation}

Next, we consider any $\rho$-pseudo-orbit $x'_k$ obtained from $x$
by adding on an error of size $\leq\rho$ at each application of the
map, i.e.~$d(x'_0,x_0)\leq \rho$, and for $k\geq 1$, $d(x'_k,Fx'_{
k-1})\leq \rho$.  Provided $d(x_j,x'_j)<\gamma $ for each $j<k$, it
follows that
\begin {equation}\label{eq:gron1}
    d(x_k,x'_k)\leq
    \rho\sum_{j=0}^{k}\left(\frac{\C}{\gamma}\right)^j\leq
    \C\,\rho \left (\frac {\C} {\gamma}\right) ^k .
\end {equation}
In particular, if $\rho $, $\gamma $, and $\lambda$ were chosen such
that
\begin {equation}\label{eq:gron2}
    \C\,\rho \left (\frac {\C} {\gamma}\right) ^\lambda<\gamma,
\end {equation}
then Equation~\eqref{eq:gron1} will hold for each $k\leq\lambda $.
We assume that Equation~\eqref{eq:gron2} is true.  Then we can also
control the differences in elapsed flight times using
Equation~\eqref{eq:2d_time_derivative}:
\begin {equation}\label{eq:gron3}
    \abs{\zeta x_k-\zeta x'_k}
    \leq
    \frac {\C\,\rho} {\gamma } \left (\frac {\C} {\gamma}\right) ^k.
\end {equation}

It remains to estimate the size $\nu \mathcal{C}_{\gamma,\lambda}$
of the set of $x$ for which the above estimates do not hold.  Using
Lemma~\ref{lem:gron1} below,
\begin {equation}\label{eq:gron4}
    \nu\mathcal{C}_{\gamma,\lambda}
    \leq
    (\lambda +1)\bigl(\nu \mathcal{N}_\gamma (\partial\Omega)+\nu \mathcal{N}_\gamma
    (F^ {-1}\mathcal{N}_\gamma
    (\partial\Omega))  \bigr)
    \leq
    \mathcal{O} (\lambda (\gamma +\gamma  ^ {1/3})) =\mathcal{O} (\lambda \gamma  ^
    {1/3}).
\end {equation}

\begin{lem}\label{lem:gron1}
    As $\gamma \rightarrow 0$,
    \[
        \nu \mathcal{N}_\gamma  (F^ {-1}\mathcal{N}_\gamma
        (\partial\Omega))=\mathcal{O}(\gamma  ^ {1/3}).
    \]
\end {lem}

This estimate is not necessarily the best possible.  For example,
for dispersing billiard tables, where the curvature of the boundary
is positive, one can show that $\nu \mathcal{N}_\gamma  (F^
{-1}\mathcal{N}_\gamma (\partial\Omega))=\mathcal{O}(\gamma  )$.
However, the estimate in Lemma~\ref{lem:gron1} is general and
sufficient for our needs.

\begin {proof}
First, we note that it is equivalent to estimate $\nu
\mathcal{N}_\gamma (F\mathcal{N}_\gamma (\partial\Omega))$, as $F$
has the measure-preserving involution $\mathcal{I} (r,\varphi) = (r,
-\varphi) $, i.e.~$F^ {-1} =\mathcal{I}\circ F\circ\mathcal{I}
$~\cite{CherMark06}.

Fix $\alpha\in (0,1/2) $, and cover $\mathcal{N}_\gamma
(\partial\Omega) $ with $\mathcal{O} (\gamma ^ {-1}) $ starlike
sets, each of diameter no greater than $\mathcal{O} (\gamma) $.  For
example, these sets could be squares of side length $\gamma $.
Enumerate the sets as $\set {A_i} $.  Set $\mathcal{G}=\set
{i:FA_i\cap \mathcal{N}_{\gamma ^\alpha} (\partial\Omega)
=\varnothing}$.

If $i\in \mathcal{G}$, $F\arrowvert _{A_i} $ is a diffeomorphism
satisfying $\norm{DF\arrowvert_{A_i}}\leq\mathcal{O} (\gamma ^
{-\alpha}) $.  See Equation~\eqref{eq:2d_derivative_bound}.  Thus
$\dia{FA_i}\leq\mathcal{O} (\gamma ^ {1-\alpha}) $, and so
$\dia{\mathcal{N}_\gamma (FA_i)}\leq\mathcal{O} (\gamma ^
{1-\alpha}) $.  Hence $\nu\mathcal{N}_\gamma (FA_i)\leq\mathcal{O}
(\gamma ^ {2(1-\alpha)}) $, and $\nu\mathcal{N}_\gamma
(\cup_{i\in\mathcal{G}}FA_i)\leq\mathcal{O} (\gamma ^ {1-2\alpha})
$.

If $i\notin \mathcal{G}$,  $A_i\cap F^ {-1} (\mathcal{N}_{\gamma
^\alpha} (\partial\Omega)) \neq\varnothing $.  Thus $A_i$ might be
cut into many pieces by $F^ {-1} (\partial\Omega) $, but each of
these pieces must be mapped near $\partial\Omega $.  In fact,
$FA_i\subset\mathcal{N}_{\mathcal{O} (\gamma ^\alpha)}
(\partial\Omega) $.  This is because outside $F^ {-1}
(\mathcal{N}_{\gamma ^\alpha} (\partial\Omega))$,
$\norm{DF}\leq\mathcal{O} (\gamma ^ {-\alpha}) $, and so points in
$FA_i$ are no more than a distance $\mathcal{O} (\gamma /\gamma ^
{\alpha}) $ away from $\mathcal{N}_{\gamma^\alpha} (\partial\Omega)
$, and $\gamma <\gamma  ^ {1-\alpha} <\gamma  ^\alpha $.  It follows
that $\mathcal{N}_\gamma (FA_i) \subset\mathcal{N}_{\mathcal{O}
(\gamma ^\alpha)} (\partial\Omega)$, and
$\nu\mathcal{N}_{\mathcal{O} (\gamma ^\alpha)}
(\partial\Omega)=\mathcal{O} (\gamma ^\alpha).  $

Thus $\nu \mathcal{N}_\gamma  (F^ {-1}\mathcal{N}_\gamma
(\partial\Omega))=\mathcal{O}(\gamma  ^ {1-2\alpha}+\gamma  ^
{\alpha})$, and we obtain the lemma by taking $\alpha=1/3$.

\end {proof}

\subsubsection{Application to a perturbed billiard flow}
\label{sct:gron_ap}

Returning to the end of Step 5 in Section~\ref{sct:main_steps}, let
the initial conditions of the slow variables be fixed at
$h_0=(Q_0,W_0,E_{1,0},E_{2,0})\in\mathcal{V} $ throughout the
remainder of this section.  We can assume that the billiard dynamics
of the left gas particle in $\mathcal{D}_1(Q_0) $ are ergodic. Also,
fix a particular value of the initial conditions for the right gas
particle for the remainder of this section.  Then $z_\varepsilon(t)
$ and $\tilde T_\varepsilon$ may be thought of as random variables
depending on the left gas particle's initial conditions
$y\in\mathcal{M} ^1$. Now if $h_\varepsilon (t)=
(Q_\varepsilon(t),W_\varepsilon(t),E_{1,\varepsilon}(t),E_{2,\varepsilon}(t))$
denotes the actual motions of the slow variables when
$\varepsilon>0$, it follows from Equation~\eqref{eq:h_div} that,
provided $\varepsilon L(\varepsilon)\leq \tilde {T}_\varepsilon $,
\begin {equation}\label{eq:h_div2}
    \sup_{0\leq t\leq
    L(\varepsilon)}\abs{h_0-h_\varepsilon(t)}=\mathcal{O}(\varepsilon L(\varepsilon)).
\end {equation}
Furthermore, we only need to show that
\begin {equation}\label{eq:gron5}
    \mu
    \set{y\in\mathcal{M} ^1:\abs{\frac{1}{L(\varepsilon)}\int_0^{L(\varepsilon)}
    G(z_\varepsilon(s))-G(z_0(s))ds}\geq\frac{\delta}{2T}
    \text { and } \varepsilon
    L(\varepsilon)\leq\tilde{T}_\varepsilon }
    \rightarrow 0
\end {equation}
as $\varepsilon\rightarrow 0$, where $G$ is defined in
Equation~\eqref{eq:G_definition}.

For definiteness, we take the following quantities from
Subsection~\ref{sct:Gronwall_map} to depend on $\varepsilon$ as
follows:
\begin {equation}\label{eq:gron6}
\begin {split}
    L(\varepsilon) &= L=\log \log\frac{1}{\varepsilon},
    \\
    \gamma (\varepsilon) &= \gamma =e^{-L},
    \\
    \lambda(\varepsilon)&=\lambda=
    \frac{2}{E_{\nu}\zeta}L,
    \\
    \rho(\varepsilon) &=\rho=\C\frac {\varepsilon L} {\gamma }.
\end {split}
\end {equation}
The constant in the choice of $\rho$ and $\rho$'s dependence on
$\varepsilon$ will be explained in the proof of
Lemma~\ref{lem:gron3}, which is at the end of this subsection. The
other choices may be explained as follows. We wish to use continuity
estimates for the billiard map to produce continuity estimates for
the flow on the time scale $L$. As the divergence of orbits should
be exponentially fast, we choose $L$ to grow sublogarithmically in
$\varepsilon^ {-1} $.  Since from Equation~\eqref{eq:2d_Santalo} the
expected flight time between collisions with
$\partial\mathcal{D}_1(Q_0)$ when $\varepsilon=0$ is
$E_{\nu}\zeta=\pi\abs{\mathcal{D}_1(Q_0)}/(\sqrt{2E_{1,0}}\abs{\partial\mathcal{D}_1(Q_0)})$,
we expect to see roughly $\lambda/2$ collisions on this time scale.
Considering $\lambda$ collisions gives us some margin for error.
Furthermore, we will want orbits to keep a certain distance, $\gamma
$, away from the billiard discontinuities. $\gamma \rightarrow 0$ as
$\varepsilon\rightarrow 0$, but $\gamma  $ is very large compared to
the possible drift $\mathcal{O} (\varepsilon L) $ of the slow
variables on the time scale $L$.  In fact, for each $C,m,n>0$,
\begin {equation}\label{eq:gron7}
    \frac {\varepsilon L^m } {\gamma ^n}
    \left (\frac {C} {\gamma}\right)^\lambda=\mathcal{O}
    (\varepsilon\, e^{\C\,L^2})
    \rightarrow 0\text { as }\varepsilon\rightarrow 0.
\end {equation}

Let $X:\mathcal{M} ^1\rightarrow\Omega$ be the map taking
$y\in\mathcal{M} ^1$ to $x=X(y)\in\Omega $, the location of the
billiard orbit of $y$ in the collision cross-section that
corresponds to the most recent time in the past that the orbit was
in the collision cross-section.  We consider the set of initial
conditions
\[
    \mathcal{E}_\varepsilon=
    X^{-1}(\Omega\backslash\mathcal{C}_{\gamma ,\lambda})\bigcap
    X^{-1}
    \set {x\in\Omega:  \sum_{k=0}^\lambda \zeta (F^k x)>
    L}.
\]
Now from Equations~\eqref{eq:gron4} and~\eqref{eq:gron6},
$\nu\mathcal{C}_{\gamma ,\lambda}\rightarrow 0$ as
$\varepsilon\rightarrow 0$.  Furthermore, by the ergodicity of $F$,
$\nu\set {x\in\Omega:\sum_{k=0}^\lambda \zeta (F^k x)\leq L}=\nu\set
{x\in\Omega:\lambda^ {-1}\sum_{k=0}^\lambda \zeta (F^k x)\leq E_\nu
\zeta/2}\rightarrow 0$ as $\varepsilon\rightarrow 0$.  But because
the free flight time is bounded above, $\mu X^ {-1}\leq \C\cdot \nu
$, and so $\mu\mathcal{E}_\varepsilon\rightarrow 1$ as
$\varepsilon\rightarrow 0$.  Hence, the convergence in
Equation~\eqref{eq:gron5} and the conclusion of the proof in
Section~\ref{sct:main_steps} follow from the lemma below and
Equation~\eqref{eq:gron7}.

\begin{lem}[Analysis of deviations along good orbits]\label{lem:gron2}
As $\varepsilon\rightarrow 0$,
\[
    \sup_{y\in\mathcal{E}_\varepsilon \cap \set{\varepsilon L\leq \tilde {T}_\varepsilon
    }}
    \abs{\frac{1}{L}\int_0^{L}
    G(z_\varepsilon(s))-G(z_0(s))ds}=
    \mathcal{O} \left(\rho
    \left(\frac {\C} {\gamma} \right)^ {\lambda}\right)
    +\mathcal{O}(L^ {- 1})\rightarrow 0.
\]

\end {lem}

\begin {proof}

Fix a particular value of $y\in\mathcal{E}_\varepsilon \cap
\set{\varepsilon L\leq \tilde {T}_\varepsilon}$. For convenience,
suppose that $y=X(y) =x\in\Omega $.  Let $y_0(t) $ denote the time
evolution of the billiard coordinates for the left gas particle when
$\varepsilon=0$.  Then there is some $N\leq \lambda$ such that the
orbit $ x_k=F^k x= ( r_k,\varphi_k)$ for $0\leq k\leq N $
corresponds to all of the instances (in order) when $y_0(t) $ enters
the collision cross-section $\Omega=\Omega_{ h_0} $ corresponding to
collisions with $\partial\mathcal{D}_1(Q_0) $ for $0\leq t\leq L$.
We write $\Omega_{ h_0}$ to emphasize that in this subsection we are
only considering the collision cross-section corresponding to the
billiard dynamics in the domain $\mathcal{D}_1(Q_0) $ at the energy
level $ E_{1,0} $.  In particular, $F$ will always refer to the
return map on $\Omega_{ h_0}$.

Also, define an increasing sequence of times $ t_k$ corresponding to
the actual times $ y_0(t) $ enters the collision cross-section, i.e.
\[
\begin {split}
    t_0 &=0,\\
    t_k & = t_{k-1} +\zeta x_{k-1}\text { for } k>0.
\end {split}
\]
Then $ x_k= y_0 ( t_k) $.  Furthermore, define inductively
\[
\begin {split}
    N_1&=\inf\set{k>0: t_k\text{ corresponds to a collision with the
    piston}},\\
    N_j&=\inf\set{k>N_{j-1}: t_k\text{ corresponds to a collision with the
    piston}}.\\
\end {split}
\]

Next, let $y_\varepsilon(t) $ denote the time evolution of the
billiard coordinates for the left gas particle when $\varepsilon>0$.
We will construct a pseudo-orbit $x_{k,\varepsilon}' =
(r_{k,\varepsilon}',\varphi_{k,\varepsilon}')$ of points in
$\Omega_{ h_0}$ that essentially track the collisions (in order) of
the left gas particle with the boundary under the dynamics of
$y_\varepsilon(t) $ for $0\leq t\leq L$.

First, define an increasing sequence of times $ t_{k,\varepsilon}'$
corresponding to the actual times $ y_\varepsilon(t) $ experiences a
collision with the boundary of the gas container or the moving
piston.  Define
\[
\begin {split}
    N_{\varepsilon}'&=\sup\set{k\geq 0: t_{k,\varepsilon}'
    \leq L},\\
    N_{1,\varepsilon}'&=\inf\set{k>0: t_{k,\varepsilon}'
    \text{ corresponds to a collision with the
    piston}},\\
    N_{j,\varepsilon}'&=\inf\set{k>N_{j-1,\varepsilon}': t_{k,\varepsilon}'
    \text{ corresponds to a collision with the
    piston}}.\\
\end {split}
\]
Because $L\leq \tilde {T}_\varepsilon(y)/\varepsilon$, we know that
as long as $N_{j+1,\varepsilon}'\leq N_{\varepsilon}'$, then
$N_{j+1,\varepsilon}'-N_{j,\varepsilon}'\geq 2$.  See the discussion
in Subsection~\ref{sct:collisions}.  Then we define
$x_{k,\varepsilon}'\in\Omega_{ h_0}$ by
\[
    x_{k,\varepsilon}' =
    \begin {cases}
    y_\varepsilon (t_{k,\varepsilon}')
    \text { if }k\notin\set {N_{j,\varepsilon}'},\\
    F^ {-1}x_{k+1,\varepsilon}'
    \text { if }k\in\set {N_{j,\varepsilon}'}.
    \end {cases}
\]

\begin{lem}\label{lem:gron3}
Provided $\varepsilon$ is sufficiently small, the following hold for
each $k\in [0,N\wedge N_\varepsilon')$.  Furthermore, the requisite
smallness of $\varepsilon$ and the sizes of the constants in these
estimates may be chosen independent of the initial condition
$y\in\mathcal{E}_\varepsilon \cap \set{\varepsilon L\leq \tilde
{T}_\varepsilon}$ and of $k$:
\begin {itemize}
    \item[\emph{(a)}]
        $x_{k,\varepsilon}' $ is well defined.  In particular, if
        $k\notin\set{N_{j,\varepsilon}'} $,
        $y_\varepsilon(t_{k,\varepsilon}') $ corresponds to a
        collision point on $\partial\mathcal{D}_1( Q_0)$, and not to
        a collision point on a piece of $\partial\mathcal{D} $ to
        the right of $ Q_0$.
    \item[\emph{(b)}]
        If $ k>0$ and $k\notin\set{N_{j,\varepsilon}'} $, then
        $x_{k,\varepsilon}' =Fx_{k-1,\varepsilon}' $.
    \item[\emph{(c)}]
        If $ k>0$ and $k\in\set{N_{j,\varepsilon}'} $, then
        $d(x_{k,\varepsilon}',Fx_{k-1,\varepsilon}')\leq\rho$ and
        the $\varphi$ coordinate of $y_\varepsilon(t_{k,\varepsilon}') $
        satisfies
        $\varphi(y_\varepsilon(t_{k,\varepsilon}')) =\varphi_{k,\varepsilon}' +
        \mathcal{O} (\varepsilon).$
    \item[\emph{(d)}]
        $d(x_k,x'_{k,\varepsilon})\leq\C\,\rho (\C/\gamma) ^k$ .
    \item[\emph{(e)}]
        $k=N_{j,\varepsilon}'$ if and only if $k=N_j$.
    \item[\emph{(f)}]
        If $ k>0$, $t_{k,\varepsilon}'-t_{k-1,\varepsilon}'
        =
        t_k - t_{ k-1} +
        \mathcal{O}(\rho \left (\C/\gamma\right) ^k).$
\end {itemize}

\end{lem}

We defer the proof of Lemma~\ref{lem:gron3} until the end of this
subsection.  Assuming that $\varepsilon$ is sufficiently small for
the conclusions of Lemma~\ref{lem:gron3} to be valid, we continue
with the proof of Lemma~\ref{lem:gron2}.

Set $M=N\wedge N_\varepsilon'-1$. Note that $M\leq\lambda\sim L $.
From (f) in Lemma~\ref{lem:gron3} and Equations~\eqref{eq:gron6} and
\eqref{eq:gron7}, we see that
\[
\begin {split}
    \abs{t_M-t_{M,\varepsilon}'}
    &
    \leq \sum_{k=1}^M
    \abs{t_{k,\varepsilon}'-t_{k-1,\varepsilon}'- (t_k - t_{ k-1})}
    =
    \mathcal{O}\left(\rho \frac{\C^\lambda}{\gamma^{\lambda}}\right)
    \rightarrow 0\text{ as }\varepsilon\rightarrow 0.
\end {split}
\]
Because the flight times $t_{k,\varepsilon}'-t_{k-1,\varepsilon}'$
and $t_k - t_{ k-1}$ are uniformly bounded above, it follows from
the definitions of $N$ and $N_\varepsilon' $ that $t_M,\,
t_{M,\varepsilon}'\geq L-\C$.  But from
Subsection~\ref{sct:collisions}, the time between the collisions of
the left gas particle with the piston are uniformly bounded away
from zero. Using (c) and Equation~\eqref{eq:h_div2}, it follows that
\[
\begin {split}
    &\abs{\frac{1}{L}
    \int_0^{L}G(z_\varepsilon(s))-G(z_0(s))ds}
    \\
    &\qquad
    =\mathcal{O} (L^ {-1}) +
    \sum_{k\in \set { N_j:N_j\leq M}} \abs{\sqrt{2E_{1,0}}\,\cos \varphi_k
    -\sqrt{2E_{1,\varepsilon}(t_{k,\varepsilon}')}\,\cos
    (\varphi_{k,\varepsilon}'+\mathcal{O} (\varepsilon))}
    \\
    &\qquad
    =\mathcal{O} (L^ {-1}) +
    \sum_{k\in \set { N_j:N_j\leq M}}
    \abs{\sqrt{2E_{1,0}}\,\cos \varphi_k
    -\sqrt{2E_{1,0}}\,\cos
    \varphi_{k,\varepsilon}'
    +\mathcal{O} (\varepsilon L)}
    \\
    &\qquad
    =\mathcal{O} (L^ {-1}) +
    \mathcal{O} (\varepsilon L^2)
    +\sqrt{2E_{1,0}}\,\sum_{k\in \set { N_j:N_j\leq M}}
    \abs{\cos \varphi_k
    -\cos
    \varphi_{k,\varepsilon}'}.
\end {split}
\]
But using (d),
\[
\begin {split}
    \sum_{k\in \set { N_j:N_j\leq M}}
    \abs{\cos \varphi_k-\cos
    \varphi_{k,\varepsilon}'}
    \leq\sum_{ k=0} ^M\mathcal{O}  (\rho (\C/\gamma) ^k)
    =\mathcal{O}  (\rho (\C/\gamma) ^\lambda).
\end {split}
\]
Since $\varepsilon L^2=\mathcal{O}(\rho (\C/\gamma) ^\lambda) $,
this finishes the proof of Lemma~\ref{lem:gron2}.

\end {proof}

\begin {proof}[Proof of Lemma~\ref{lem:gron3}]

The proof is by induction.  We take $\varepsilon$ to be so small
that Equation~\eqref{eq:gron2} is satisfied. This is possible by
Equation~\eqref{eq:gron7}.

It is trivial to verify (a)-(f) for $ k=0$. So let $0<l<N\wedge
N_\varepsilon' $, and suppose that (a)-(f) have been verified for
all $k<l$.  We have three cases to consider:

\subsubsection*{Case 1: $l-1$ and $l\notin \set{N_{j,\varepsilon}'}$:}

In this case, verifying (a)-(f) for $ k=l$ is a relatively
straightforward application of the machinery developed in
Subsection~\ref{sct:Gronwall_map}, because for
$t_{l-1,\varepsilon}'\leq t\leq t_{l,\varepsilon}'$, $y_\varepsilon
(t) $ traces out the billiard orbit between $x_{l-1,\varepsilon}'$
and $x_{l,\varepsilon}'$ corresponding to free flight in the domain
$\mathcal{D}_1( Q_0) $. We make only two remarks.

First, as long as $\varepsilon$ is sufficiently small, it really is
true that $x_{l,\varepsilon}'=y_\varepsilon (t_{l,\varepsilon}')$
corresponds to a true collision point on $\partial\mathcal{D}_1(
Q_0) $.  Indeed, if this were not the case, then it must be that
$Q_\varepsilon(t_{l,\varepsilon}')> Q_0 $, and $y_\varepsilon
(t_{l,\varepsilon}')$ would have to correspond to a collision with
the side of the ``tube'' to the right of $ Q_0$.  But then
$x_{l,\varepsilon}''  =Fx_{l-1,\varepsilon}'\in\Omega_{ h_0}$ would
correspond to a collision with an immobile piston at $ Q_0$ and
would satisfy $d(x_k,x''_{k,\varepsilon})\leq\C\,\rho (\C/\gamma)
^k\leq \C\,\rho (\C/\gamma) ^\lambda =o(\gamma )$, using
Equations~\eqref{eq:gron1} and \eqref{eq:gron7}.  But $
x_k\notin\mathcal{N}_\gamma (\partial\Omega_{ h_0}) $, and so it
follows that when the trajectory of $y_\varepsilon(t) $ crosses the
plane $\set {Q= Q_0} $, it is at least a distance $\sim \gamma $
away from the boundary of the face of the piston, and its velocity
vector is pointed no closer than $\sim \gamma $ to being parallel to
the piston's face.  As $Q_\varepsilon(t_{l,\varepsilon}')-
Q_0=\mathcal{O} (\varepsilon L) =o(\gamma ) $, and it is
geometrically impossible (for small $\varepsilon$) to construct a
right triangle whose sides $ s_1,\: s_2$ satisfy $\abs{ s_1}\geq\sim
\gamma ,\:\abs{ s_2}\leq\mathcal{O} (\varepsilon L) $, with the
measure of the acute angle adjacent to $ s_1$ being greater than
$\sim \gamma$, we have a contradiction. After crossing the plane
$\set {Q= Q_0} $, $y_\varepsilon(t) $ must experience its next
collision with the face of the piston, which violates the fact that
$l\notin \set{N_{j,\varepsilon}'}$.

Second, $t_{l,\varepsilon}'-t_{l-1,\varepsilon}' =\zeta
x'_{l-1,\varepsilon}+\mathcal{O}(\varepsilon L)$, because
$v_{1,\varepsilon} =v_ {1,0} +\mathcal{O} (\varepsilon L) $.  See
Equation~\eqref{eq:h_div2}.  From Equation~\ref{eq:gron3},
$\abs{\zeta x_{l-1}-\zeta x'_{l-1,\varepsilon}} \leq \mathcal{O}
((\rho/\gamma)  \left (\C/\gamma\right) ^{l-1})$. As $t_l - t_{
l-1}=\zeta x_{l-1}$ and $\varepsilon L=\mathcal{O} ((\rho/\gamma)
\left (\C/\gamma\right) ^{l-1})$, we obtain (f).

\subsubsection*{Case 2: There exists $i$ such that $l=N_{i,\varepsilon}'$: }

For definiteness, we suppose that
$Q_\varepsilon(t_{l,\varepsilon}')\geq Q_0$, so that the left gas
particle collides with the piston to the right of $ Q_0$.  The case
when $Q_\varepsilon(t_{l,\varepsilon}')\leq Q_0$ can be handled
similarly.

We know that $ x_{l-1}, x_{l},x_{l+1}\notin \mathcal{N}_\gamma
(\partial\Omega_{ h_0})\cup\mathcal{N}_\gamma  (F^
{-1}\mathcal{N}_\gamma (\partial\Omega_{ h_0}))$.  Using the
inductive hypothesis and Equation~\eqref{eq:gron1}, we can define
\[
    x_{l,\varepsilon}''=Fx_{l-1,\varepsilon}',\qquad
    x_{l+1,\varepsilon}''=F^2x_{l-1,\varepsilon}',
\]
and $d( x_{l},x_{l,\varepsilon}'') \leq\C\,\rho (\C/\gamma) ^{l}$,
$d( x_{l+1}, x_{l+1,\varepsilon}'')\leq\C\,\rho (\C/\gamma) ^{l+1}$.
In particular, $x_{l,\varepsilon}''$ and $x_{l+1,\varepsilon}''$ are
both a distance $\sim \gamma $ away from $\partial\Omega_{ h_0} $.
Furthermore, when the left gas particle collides with the moving
piston, it follows from Equation~\eqref{eq:v_1Wchange} that the
difference between its angle of incidence and its angle of
reflection is $\mathcal{O} (\varepsilon) $.  Referring to
Figure~\ref{fig:collision}, this means that
$\varphi_{l,\varepsilon}' =\varphi_{l,\varepsilon}'' +\mathcal{O}
(\varepsilon) $. Geometric arguments similar to the one given in
Case 1 above show that the $y_\varepsilon$-trajectory of the left
gas particle has precisely one collision with the piston and no
other collisions with the sides of the gas container when the gas
particle traverses the region $Q_0\leq Q\leq
Q_\varepsilon(t_{l,\varepsilon}')$. Note that $x_{l,\varepsilon}' $
was defined to be the point in the collision cross-section $\Omega_{
h_0} $ corresponding to the return of the $y_\varepsilon$-trajectory
into the region $Q\leq Q_0$. See Figure~\ref{fig:collision}.  From
this figure, it is also evident that
$d(r_{l,\varepsilon}',r_{l,\varepsilon}'')\leq\mathcal{O}
(\varepsilon L/\gamma ) $.  Thus $d(
x_{l,\varepsilon}'',x_{l,\varepsilon}')=\mathcal{O} (\varepsilon
L/\gamma )$, and this explains the choice of $\rho(\varepsilon) $ in
Equation~\eqref{eq:gron6}.

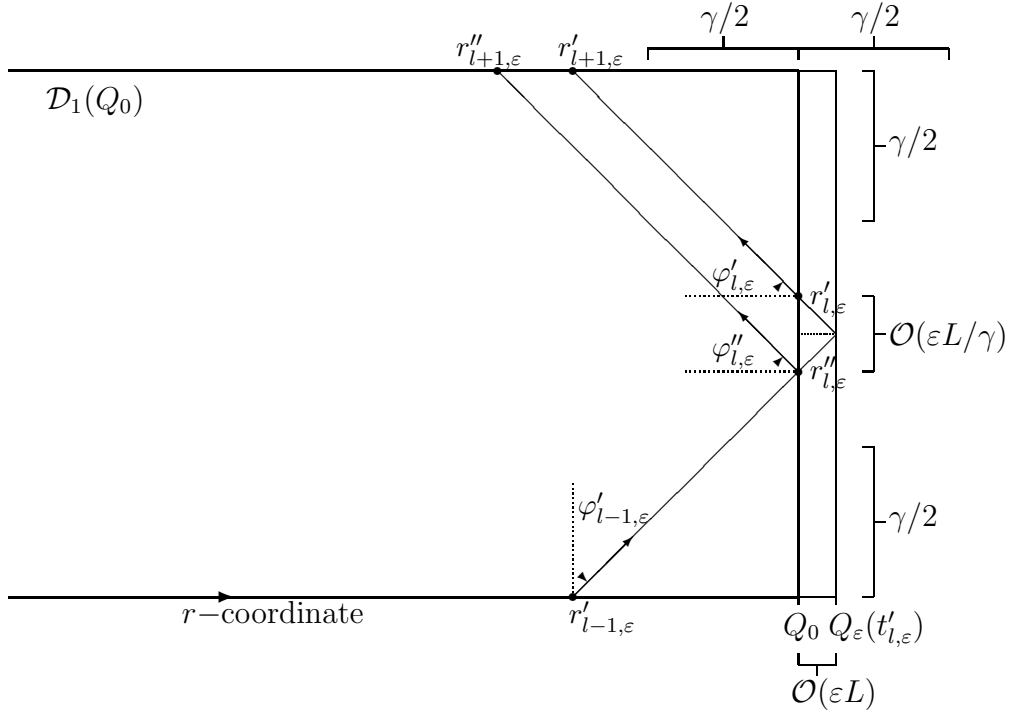
\begin{figure}
    \begin {center}
    \setlength{\unitlength}{1.0 cm}
    \begin{picture}(15,10)
        \thicklines
        \put(1,1.5){\line(1,0){10.5}}
        \put(1,8.5){\line(1,0){10.5}}
        \put(11.5,1.5){\line(0,1){7.0}}
        \put(3.5,1.5){\vector(1,0){0.5}}
        \put(3.3,1.15){$r-$coordinate}
        \thinlines
        \put(1,1.5){\line(1,0){11}}
        \put(1,8.5){\line(1,0){11}}
        \put(11.5,1.4){\line(0,1){7.1}}
        \put(11.3,1.0){$Q_0$}
        \put(12,1.4){\line(0,1){7.1}}
        \put(11.9,1.0){$Q_\varepsilon(t_{l,\varepsilon}')$}
        \put(1.5,8){$\mathcal{D}_1(Q_0) $}
        \put(9.5,8.8){\line(1,0){4}}
        \put(11.5,8.8){\line(0,-1){0.15}}
        \put(13.5,8.8){\line(0,-1){0.15}}
        \put(9.5,8.8){\line(0,-1){0.15}}
        \put(  10.5,8.8){\line(0,1){0.15}}
        \put(10.2,9.1){$\gamma/2 $}
        \put(12.5,8.8){\line(0,1){0.15}}
        \put(12.2,9.1){$\gamma/2 $}
        \put(11.5,0.6){\line(1,0){0.5}}
        \put(11.75,0.6){\line(0,-1){0.15}}
        \put(  11.5,0.6){\line(0,1){0.15}}
        \put(  12,0.6){\line(0,1){0.15}}
        \put(11.4,0.1){$\mathcal{O} (\varepsilon L) $}
        \put(12.5,4.5){\line(0,1){1}}
        \put(12.5,4.5){\line(-1,0){0.15}}
        \put(12.5,5){\line(1,0){0.15}}
        \put(12.5,5.5){\line(-1,0){0.15}}
        \put(12.7,4.85){$\mathcal{O} (\varepsilon L/\gamma ) $}
        \put(12.5,6.5){\line(0,1){2}}
        \put(12.5,8.5){\line(-1,0){0.15}}
        \put(12.5,7.5){\line(1,0){0.15}}
        \put(12.5,6.5){\line(-1,0){0.15}}
        \put(12.7,7.45){$\gamma /2 $}
        \put(12.5,1.5){\line(0,1){2}}
        \put(12.5,3.5){\line(-1,0){0.15}}
        \put(12.5,2.5){\line(1,0){0.15}}
        \put(12.5,1.5){\line(-1,0){0.15}}
        \put(12.7,2.45){$\gamma  /2$}
    \small
        \put(12,5){\line(-1,-1){3.5}}
        \put(11.5,4.5){\vector(-1,1){0.8}}
        \put(11.5,4.5){\line(-1,1){4}}
        \put(11.65,4.4){$r_{l,\varepsilon}''$}
        \put(10.35,4.7){$\varphi_{l,\varepsilon}''$}
        \put(11.5,4.5){\circle*{.1}}
        \qbezier[24](10,4.5)(10.75,4.5)(11.5,4.5)
        \put(11.1,4.5){\vector(1,1){0.2}}
        \put(12,5){\vector(-1,1){1.3}}
        \put(12,5){\line(-1,1){3.5}}
        \put(11.65,5.4){$r_{l,\varepsilon}'$}
        \put(10.35,5.7){$\varphi_{l,\varepsilon}'$}
        \put(11.5,5.5){\circle*{.1}}
        \qbezier[24](10,5.5)(10.75,5.5)(11.5,5.5)
        \put(11.1,5.5){\vector(1,1){0.2}}
        \qbezier[12](11.5,5)(11.75,5)(12,5)
        \put(8.5,1.5){\circle*{.1}}
        \put(8.5,1.5){\vector(1,1){0.8}}
        \qbezier[24](8.5,1.5)(8.5,2.25)(8.5,3.0)
        \put(8.45,1.15){$r_{l-1,\varepsilon}'$}
        \put(8.57,2.6){$\varphi_{l-1,\varepsilon}'$}
        \put(8.5,1.9){\vector(1,-1){0.2}}
        \put(8.5,8.5){\circle*{.1}}
        \put(8.3,8.7){$r_{l+1,\varepsilon}'$}
        \put(7.5,8.5){\circle*{.1}}
        \put(6.95,8.7){$r_{l+1,\varepsilon}''$}
    \normalsize
    \end{picture}
    \end {center}
    \caption{An analysis of the divergences of orbits when $\varepsilon>0$
        and the left gas particle collides with the moving piston to the right of $Q_0$.  Note that the
        dimensions are distorted for visual clarity, but that $\varepsilon L$
        and $\varepsilon L/\gamma $ are both $o(\gamma ) $ as $\varepsilon\rightarrow
        0$.}  Furthermore, $\varphi_{l,\varepsilon}''\in(-\pi/2+\gamma/2,\pi/2-\gamma/2) $
        and $\varphi_{l,\varepsilon}' =\varphi_{l,\varepsilon}'' +\mathcal{O} (\varepsilon) $,
        and so $r_{l,\varepsilon}' =r_{l,\varepsilon}''+\mathcal{O}
        (\varepsilon L/\gamma ) $.  In particular, the
        $y_\varepsilon$-trajectory of the left gas particle has precisely
        one collision with the piston and no other collisions with the sides
        of the gas container when the gas particle traverses the region
        $Q_0\leq Q\leq Q_\varepsilon(t_{l,\varepsilon}')$
    \label{fig:collision}
\end{figure}

From the above discussion and the machinery of
Subsection~\ref{sct:Gronwall_map}, (a)-(e) now follow readily for
\emph{both} $ k=l$ and $ k=l+1$.  Furthermore, property (f) follows
in much the same manner as it did in Case 1 above.  However, one
should note that $t_{l,\varepsilon}'-t_{l-1,\varepsilon}' =\zeta
x'_{l-1,\varepsilon}+\mathcal{O}(\varepsilon L
)+\mathcal{O}(\varepsilon L/\gamma )$ and
$t_{l+1,\varepsilon}'-t_{l,\varepsilon}' =\zeta
x'_{l,\varepsilon}+\mathcal{O}(\varepsilon L
)+\mathcal{O}(\varepsilon L/\gamma )$, because of the extra distance
$\mathcal{O} (\varepsilon L/\gamma ) $ that the gas particle travels
to the right of $ Q_0$.  But $\varepsilon L/\gamma =\mathcal{O}
((\rho/\gamma) \left (\C/\gamma\right) ^{l-1})$, and so property (f)
follows.

\subsubsection*{Case 3: There exists $i$ such that $l-1=N_{i,\varepsilon}'$: }

As mentioned above, the inductive step in this case follows
immediately from our analysis in Case 2.

\end {proof}

\section{Generalization to a full proof of Theorem~\ref{thm:dDpiston}}
\label{sct:generalization}

It remains to generalize the proof in Sections~\ref{sct:2dprep} and
\ref{sct:2dproof} to the cases when $ n_1, n_2\geq 1$ and $d=3$.

\subsection{Multiple gas particles on each side of the piston}\label{sct:multiple}

When $d=2$, but $n_1,n_2\geq 1$, only minor modifications are
necessary to generalize the proof above. As in
Subsection~\ref{sct:collisions}, one defines a stopping time $\tilde
{T}_\varepsilon$ satisfying $P\set {\tilde {T}_\varepsilon < T\wedge
T_\varepsilon} =\mathcal{O} (\varepsilon) $ such that for $0\leq
t\leq \tilde {T}_\varepsilon/\varepsilon$, gas particles will only
experience clean collisions with the piston.

Next, define $H(z) $ by
\[
    H(z) =
    \begin{bmatrix}
    W\\
    +2\sum_{j=1}^{n_1}\abs{v_{1,j}^{\perp }} \delta_{q_{1,j}^{\perp }=Q}
    -2\sum_{j=1}^{n_2}\abs{v_{2,j}^{\perp }} \delta_{q_{2,j}^{\perp }=Q}\\
    -2W\abs{v_{1,j}^{\perp }} \delta_{q_{1,j}^{\perp }=Q}\\
    +2W\abs{v_{2,j}^{\perp }} \delta_{q_{2,j}^{\perp }=Q}\\
    \end{bmatrix}.
\]
It follows that for $0 \leq t\leq \tilde
{T}_\varepsilon/\varepsilon$, $
    h_\varepsilon(t)-h_\varepsilon(0)=
    \mathcal{O}(\varepsilon)+\varepsilon\int_0^t
    H(z_\varepsilon(s))ds.
$ From here, the rest of the proof follows the same steps made in
Subsection~\ref{sct:main_steps}.  We note that at Step 3, we find
that $H(z) -\bar H(h(z)) $ divides into $n_1+ n_2 $ pieces, each of
which depends on only one gas particle when the piston is held
fixed.

\subsection{Three dimensions}\label{sct:higher_d}

The proof of Theorem~\ref{thm:dDpiston} in $d=3$ dimensions is
essentially the same as the proof in two dimensions given above. The
principal differences are due to differences in the geometry of
billiards.  We indicate the necessary modifications.

In analogy with Section~\ref{sct:billiard}, we briefly summarize the
necessary facts for the billiard flows of the gas particles when
$M=\infty $ and the slow variables are held fixed at a specific
value $h\in\mathcal{V} $.  As before, we will only consider the
motions of one gas particle moving in $\mathcal{D}_1 $.  Thus we
consider the billiard flow of a point particle moving inside the
domain $\mathcal{D}_1$ at a constant speed $\sqrt{2E_1} $.  Unless
otherwise noted, we use the notation from
Section~\ref{sct:billiard}.

The billiard flow takes place in the five-dimensional space
$\mathcal{M}^1=\{(q_1,v_1)\in\mathcal{TD}_1:q_1\in\mathcal{ D}_1,\;
\abs{v_1}=\sqrt{2E_1}\}/\sim$.  Here the quotient means that when
$q_1\in\partial\mathcal{ D}_1$, we identify velocity vectors
pointing outside of $\mathcal{D}_1$ with those pointing inside
$\mathcal{D}_1$ by reflecting orthogonally through the tangent plane
to $\partial\mathcal{D}_1$ at $ q_1$. The billiard flow preserves
Liouville measure restricted to the energy surface.  This measure
has the density $d\mu=dq_1dv_1/(8\pi E_1\abs{\mathcal{D}_1} ) $.
Here $dq_1$ represents volume on $\mathbb{R}^3$, and $ dv_1$
represents area on
$S^2_{\sqrt{2E_1}}=\set{v_1\in\mathbb{R}^3:\abs{v_1}=\sqrt{2E_1}}$.

The collision cross-section
$\Omega=\{(q_1,v_1)\in\mathcal{TD}_1:q_1\in\partial\mathcal{ D}_1,\;
\abs{v_1}=\sqrt{2E_1}\}/\sim$ is properly thought of as a fiber
bundle, whose base consists of the smooth pieces of
$\partial\mathcal{D}_1$ and whose fibers are the set of outgoing
velocity vectors at $q_1\in\partial\mathcal{ D}_1$.  This and other
facts about higher-dimensional billiards, with emphasis on the
dispersing case, can be found in~\cite{BalCheSzaTot_2003}.  For our
purposes, $\Omega$ can be parameterized as follows. We decompose
$\partial\mathcal{D}_1$ into a finite union $\cup_j \Gamma_j$ of
pieces, each of which is diffeomorphic via coordinates $r$ to a
compact, connected subset of $\mathbb{R}^2$ with a piecewise
$\mathcal{C} ^3$ boundary.  The $\Gamma_j$ are nonoverlapping,
except possibly on their boundaries. Next, if $(q_1,v_1)\in\Omega $
and $ v_1$ is the outward going velocity vector, let $\hat v =
v_1/\abs{v_1} $.  Then $\Omega$ can be parameterized by $\{x=(r,\hat
v)\}$. It follows that $\Omega$ it is diffeomorphic to $\cup_j
\Gamma_j\times S^{2 +}$, where $S^{2 +}$ is the upper unit
hemisphere, and by $\partial\Omega$ we mean the subset diffeomorphic
to $(\cup_j
\partial\Gamma_j\times S^{2 +})\bigcup (\cup_j \Gamma_j\times
\partial S^{2 +})$. If $x\in\Omega $, we let $\varphi\in [0,\pi /2]$
represent the angle between the outgoing velocity vector and the
inward pointing normal vector $n$ to $\partial\mathcal{D}_1$,
i.e.~$\cos\varphi=\langle \hat v, n\rangle$. Note that we no longer
allow $\varphi$ to take on negative values. The return map
$F:\Omega\circlearrowleft$ preserves the projected probability
measure $\nu $, which has the density $d\nu=\cos\varphi\, d\hat v \,
dr/(\pi\abs{\partial\mathcal{D}_1}) $. Here
$\abs{\partial\mathcal{D}_1}$ is the area of
$\partial\mathcal{D}_1$.

$F$ is an invertible, measure preserving transformation that is
piecewise $\mathcal{C} ^2$.  Because of our assumptions on
$\mathcal{D}_1$, the free flight times and the curvature of
$\partial\mathcal{D}_1$ are uniformly bounded.  The bound on $\norm
{DF(x)}$ given in Equation~\eqref{eq:2d_derivative_bound} is still
true.  A proof of this fact for general three-dimensional billiard
tables with finite horizon does not seem to have made it into the
literature, although see~\cite{BalCheSzaTot_2003} for the case of
dispersing billiards.  For completeness, we provide a sketch of a
proof for general billiard tables in Appendix~\ref{sct:d_bounds}.

We suppose that the billiard flow is ergodic, so that $F$ is
ergodic.  Again, we induce $F$ on the subspace $\hat\Omega$ of
$\Omega$ corresponding to collisions with the (immobile) piston to
obtain the induced map $\hat F:\hat\Omega\circlearrowleft$ that
preserves the induced measure $\hat \nu$.

The free flight time $\zeta:\Omega\rightarrow \mathbb{R}$ again
satisfies the derivative bound given in
Equation~\eqref{eq:2d_time_derivative}. The generalized
Santal\'{o}'s formula\cite{Chernov1997} yields
\[
    E_\nu \zeta=\frac {4
    \abs{\mathcal{D}_1}} {\abs{v_1}\abs{\partial\mathcal{D}_1}}.
\]
If $\hat\zeta:\hat\Omega\rightarrow\mathbb{R} $ is the free flight
time between collisions with the piston, then it follows from
Proposition \ref{prop:inducing} that
\[
    E_{\hat\nu} \hat\zeta=\frac {4
    \abs{\mathcal{D}_1}}{\abs{v_1}\ell}.
\]

The expected value of $ \abs{v_1^\perp }$ when the left gas particle
collides with the (immobile) piston is given by
\[
    E_{\hat\nu} \abs{v_1^\perp }=E_{\hat\nu} \sqrt{2E_1}\cos\varphi=
    \frac{\sqrt{2E_1}}{\pi}\iint_{S^ {2+}} \cos^2\varphi\,d\hat v_1=
    \sqrt{2E_1}\frac{2}{3}.
\]

As a consequence, we obtain
\begin {lem}
\label{lem:ae_convergence_3d}

For $\mu-a.e.$ $y\in \mathcal{M}^1$,
\[
    \lim_{t\rightarrow\infty} \frac{1}{t}
    \int_0^t \abs{v_1^\perp (s)}\delta_{q_1^\perp (s)
    =Q}ds=
    \frac{E_1\ell}{3\abs{\mathcal{D}_1(Q)}}.
\]

\end {lem}
\noindent Compare the proof of Lemma~\ref{lem:ae_convergence}.

With these differences in mind, the rest of the proof of
Theorem~\ref{thm:dDpiston} when $d=3$ proceeds in the same manner as
indicated in Sections~\ref{sct:2dprep}, \ref{sct:2dproof} and
\ref{sct:multiple} above.  The only notable difference occurs in the
proof of the Gronwall-type inequality for billiards.  Due to
dimensional considerations, if one follows the proof of
Lemma~\ref{lem:gron1} for a three-dimensional billiard table, one
finds that $\nu \mathcal{N}_\gamma  (F^ {-1}\mathcal{N}_\gamma
(\partial\Omega))=\mathcal{O}(\gamma  ^ {1-4\alpha}+\gamma  ^
{\alpha})$.  The optimal value of $\alpha$ is $1/5$, and so $ \nu
\mathcal{N}_\gamma  (F^ {-1}\mathcal{N}_\gamma
(\partial\Omega))=\mathcal{O}(\gamma  ^ {1/5})$ as $\gamma
\rightarrow 0$.  Hence $\nu\mathcal{C}_{\gamma,\lambda} =\mathcal{O}
(\lambda \gamma  ^ {1/5})$, which is a slightly worse estimate than
the one in Equation~\eqref{eq:gron4}.  However, it is still
sufficient for all of the arguments in Section~\ref{sct:gron_ap},
and this finishes the proof.

\comment{
 can also be heuristically justified by
the procedure we used to justify the averaged equation when $d=2$.
Let $B^d=\{(x_1,\dots,x_d)\in\mathbb{R}^d :\sum_{i=1}^{d}x_i^2\leq
1\}$ denote the unit ball in $\mathbb{R}^d$, and let
$S^{d-1}=\{(x_1,\dots,x_d)\in B^d :\sum_{i=1}^{d}x_i^2= 1\}$ denote
the unit $(d-1)$-sphere.  Also let $(S^{d-1}) ^
+=\{(x_1,\dots,x_d)\in S^{d-1} :x_d\geq 0\}$.  Then if the piston is
held fixed, the expected flight time between collisions for the left
gas particle (with respect to the invariant billiard measure for the
billiard map induced on the subspace of collisions with the piston)
is
\[
    E_{\hat\nu} \hat\zeta=\frac {
    \abs{\mathcal{D}_1}}{\ell\sqrt{2E_1}}\frac {\abs{S^{d-1}}}
    {\abs{B^{d-1}}}.
\]
See \cite{CM06}. Furthermore,\marginal {I should say a word about
invariant measures/how to derive these equations} For future
reference, we observe that the expected value of $ \abs{v_1^\perp }$
when the left gas particle collides with the (immobile) piston is
given by
\[
    E_{\hat\nu} \abs{v_1^\perp }=
    \frac {\sqrt{2E_1}} {2d}\frac{\abs{S^{d-1}}}{\abs{B^{d-1}}}.
\]
} 

\appendix

\section{Inducing maps on subspaces}
\label{sct:inducing}

Here we present some well-known facts on inducing measure preserving
transformations on subspaces. Let $F: (\Omega,
\mathfrak{B},\nu)\circlearrowleft$ be an invertible, ergodic,
measure preserving transformation of the probability space $\Omega$
endowed with the $\sigma$-algebra $\mathfrak{B}$ and the probability
measure $\nu$.  Let $\hat\Omega\in\mathfrak{B}$ satisfy $0<\nu
\hat\Omega<1$.  Define $R:\hat\Omega\rightarrow\mathbb{N}$ to be the
first return time to $\hat\Omega$, i.e.~$R\omega
=\inf\{n\in\mathbb{N}:F^n\omega \in\hat\Omega\}$.  Then if
$\hat{\nu} : =\nu(\cdot\cap\hat\Omega)/\nu\hat\Omega$ and
$\hat{\mathfrak{B}}: =\{B\cap\hat\Omega:B\in \mathfrak{B}\}$,
$\hat{F}: (\hat\Omega,
\hat{\mathfrak{B}},\hat{\nu})\circlearrowleft$ defined by
$\hat{F}\omega=F^{R\omega}\omega$ is also an invertible, ergodic,
measure preserving transformation~\cite{Pet83}. Furthermore
$E_{\hat{\nu}} R=\int_{\hat\Omega}
R\,d\hat{\nu}=(\nu\hat\Omega)^{-1}$.

This last fact is a consequence of the following proposition:
\begin {prop}
\label{prop:inducing}

If $\zeta:\Omega\rightarrow\mathbb{R}_{\geq 0}$ is in $L^1(\nu)$,
then $\hat\zeta =\sum_{n=0}^{R-1}\zeta\circ F^n$ is in
$L^1(\hat{\nu})$, and
\[
    E_{\hat{\nu}} \hat\zeta =\frac {1}{\nu\hat\Omega}
    E_{\nu}\zeta.
\]
\end{prop}

\begin {proof}
\[
\begin {split}
    \nu\hat\Omega \int_{\hat\Omega} \sum_{n=0}^{R-1}\zeta\circ F^n\,d\hat{\nu}
    &=
    \int_{\hat\Omega} \sum_{n=0}^{R-1}\zeta\circ F^n\,d\nu
    =
    \sum_{k=1}^\infty\int_{\hat\Omega\cap\{R=k\}}
    \sum_{n=0}^{k-1}\zeta\circ F^n\,d\nu
    \\
    &=
    \sum_{k=1}^\infty\sum_{n=0}^{k-1}\int_{F^n(\hat\Omega\cap\{R=k\})}
    \zeta\,d\nu
    =
    \int_{\Omega}\zeta\,d\nu,
\end {split}
\]
because $\{F^n(\hat\Omega\cap\{R=k\}):0\leq n< k<\infty\}$ is a
partition of $\Omega$.

\end {proof}

\section{Derivative bounds for the billiard map in three dimensions}
\label{sct:d_bounds}

Returning to Section~\ref{sct:higher_d}, we need to show that for a
billiard table $\mathcal{D}_1\subset\mathbb{R}^3$ with a piecewise
$\mathcal{C} ^3$ boundary and the free flight time uniformly bounded
above, the billiard map $F$ satisfies the following: If $x_0\notin
\partial\Omega \cup F^{-1} (\partial\Omega) $, then
\begin{equation*}
    \norm {DF(x_0)}\leq\frac {\C} {\cos \varphi(Fx_0)}.
\end{equation*}

Fix $x_0= (r_0,\hat v_0)\in\Omega $, and let $x_1= (r_1,\hat
v_1)=Fx_0$.  Let $\Sigma$ be the plane that perpendicularly bisects
the straight line between $ r_0$ and $ r_1$, and let $r_{1/2} $
denote the point of intersection.  We consider $\Sigma$ as a
``transparent'' wall, so that in a neighborhood of $ x_0$, we can
write $F=F_2\circ F_1$.  Here, $F_1$ is like a billiard map in that
it takes points (i.e.~directed velocity vectors with a base) near $
x_0$ to points with a base on $\Sigma$ and a direction pointing near
$ r_1$. ($ F_1$ would be a billiard map if we reflected the image
velocity vectors orthogonally through $\Sigma$.)  $ F_2$ is a
billiard map that takes points in the image of $F_1$ and maps them
near $ x_1$. Let $x_{1/2} = F_1 x_0= F_2^ {-1} x_1 $. Then $ \norm
{DF(x_0)}\leq \norm {DF_1(x_0)}\norm {DF_2(x_{1/2})}$.

It is easy to verify that $\norm {DF_1(x_0)}\leq\C$, with the
constant depending only on the curvature of $\partial\mathcal{D}_1$
at $ r_0$.  In other words, the constant may be chosen independent
of $ x_0$. Similarly, $\norm {DF_2^ {-1}(x_1)}\leq\C$.  Because
billiard maps preserve a probability measure with a density
proportional to $\cos\varphi $, $\text {det}DF_2^ {-1}(x_1)=\cos
\varphi_{1}/\cos\varphi_{ 1/2} =\cos\varphi_1$.  As $\Omega$ is $4
$-dimensional, it follows from Cramer's Rule for the inversion of
linear transformations that
\[
    \norm {DF_2(x_{1/2})}\leq \frac {\C\norm {DF_2^ {-1}(x_1)}^3}
    {\text {det}DF_2^ {-1}(x_1)}\leq\frac {\C} {\cos\varphi_1},
\]
and we are done.

\vskip 1cm

\textbf{Acknowledgments.} The author is grateful to D. Dolgopyat,
who first introduced him to this problem, and who generously shared
his unpublished notes on averaging~\cite{Dol05}. The author also
thanks L.-S.~Young for useful discussions regarding this project and
P.~Balint for many helpful comments on the manuscript.  This
research was partially supported by the National Science Foundation
Graduate Research Fellowship Program.

\bibliographystyle{alpha}

\end{document}